\documentclass[a4paper, BCOR=0.0mm, DIV=calc]{scrartcl}
\usepackage[T1]{fontenc}
\usepackage[utf8]{inputenc}
\usepackage{url}
\usepackage{lmodern}
\usepackage{microtype}
\usepackage[english]{babel}
\usepackage{amsmath, amssymb, amsfonts}
\usepackage{nicefrac}
\newcommand{\diff}[1]{\mathrm{#1}}

\KOMAoptions{twoside=false, twocolumn=false, headinclude=false, footinclude=false, mpinclude=false, pagesize=auto}
\recalctypearea

\begin{document}
\title{Translation of C. J. Malmstèn's paper "De integralibus quibusdam definitis, seriebusque infinitis"}
\author{Alexander Aycock}
\date{}
\maketitle
\begin{abstract}
We will provide a translation of Malmstén's paper "De integralibus quibusdam definitis, seriebusque infinitis"
\end{abstract}

\section*{~~~~~~~~~~~~~~~~~~~~~~~~~~~~Introduction}

In 1846 (published in the Crelle Journal für die reine und angewandte Mathematik) the nowadays widely unknown Swedish mathematician C. J. Malmstén wrote his quite remarkable paper "De integralibus quibusdam definitis, seriebusque infinitis" \cite{9}, containing many interesting results; for example remarkable evaluations of special kinds of integrals, the de facto first proof of the Fourier series expansion for $\log{\Gamma{(x})}$ - where $\Gamma({x})$ is the well-known Gamma function, defined as $\Gamma({x}) = \int_0^{\infty}t^{x-1}e^{-t}\diff{d}t$ - (Kummer, due to whom this result is often attributed - see his paper "Beitrag zur Theorie der Function $\Gamma({x}) = \int_0^{\infty}e^{-v}v^{x-1}\diff{d}v$" \cite{7} - gave his result one year later in 1847) and maybe, being the most interesting result, gave the first proofs for functional equations of certain Dirichlet series. Most of these results are also attributed to others, as we will see later. Therefore this paper is, not only from a historical point of view, worth to be studied and  to be saved form oblivion. So I will - as I said in the paper "Note on Malmstén's paper "De integralibus quibusdam definitis, seriebusque infinitis"" \cite{11} - translate the paper from the Latin original into English as good as I can. Whenever it seems appropiate, I will make some additional notes, to give some information and/or render things more clear.\\

It is indeed quite hard to gather some information about C.J. Malmstèn, the only source, that I could locate and that describes his life and work, is " Mathematical and Mathematicians: Mathematics in Sweden Before 1950 (History of Mathematics)" by Lars Garding \cite{5}. So anyone more interested in Malmstén's background and his other works, is refered to this book, where a chapter is dedicated to him. But now we will go on to the translation. 

\section*{~~~~~On certain definite integrals and infinite series}
\subsection*{~~~written by C. J. Malmstèn, Prof. of Mathematics in Upsala}

\section*{~~~~~~~~~~~~~~~~~~~~~~~~~~~~~~~~~§1}

If we multiply the known formula
\[
\int_0^{\infty} e^{-xz}\sin{(uz)}\diff{d}z=\frac{u}{x^2+u^2}
\]
by $\diff{d}u$ and integrate from $u=0$, we will have
\[
\frac{1}{2}.~~~~2 \int_0^{\infty}\frac{e^{-xz}-e^{-xz}\cos{(uz)}}{z}\cdot\diff{d}z=\log{(x^2+u^2)}-2\log{x}
\]
whence, because it is
\[
\log{x}=\int_0^{\infty}\frac{e^{-z}-e^{-xz}}{z}\cdot\diff{d}z
\]
one will obtain
\[
1.~~~~~~ \log{(x^2+u^2)}=2 \int_0^{\infty}\frac{e^{-z}-e^{-xz}\cos{(uz)}}{z}\cdot\diff{d}z
\]
Now let us multiply on both sides by
\[
\frac{e^{au}-e^{-au}}{e^{\pi u}-e^{-\pi u}}\cdot\diff{d}u,~~~ [a<\pi]
\]
after having integrated between $u=0$ and $u=\infty$, this becomes
\[
\int_0^{\infty} \frac{e^{au}-e^{-au}}{e^{\pi u}-e^{-\pi u}} \log{(x^2+u^2)} \diff{d}u=
2 \int_0^{\infty} \frac{\diff{d}z}{z} \int_0^{\infty} \frac{e^{au}-e^{-au}}{e^{\pi u}-e^{-\pi u}}(e^{-z}-e^{-xz}\cos{(uz)}\diff{d}z
\]
whence by the means of known formula \footnote{See Exerc. de. Calc. Integra. par Legendre \cite{8}, Tom. 2, page 186 C.J. Malmstèn} if follows
\[
\int_0^{\infty} \frac{e^{au}-e^{-au}}{e^{\pi u}-e^{-\pi u}} \log{(x^2+u^2)} \diff{d}u=
\int_0^{\infty} \frac{e^{-z}\diff{d}z}{z}\bigg[\tan{\frac{1}{2}a}-\frac{2e^{-xz}\sin{a}}{1+2e^{-z}\cos{a}+e^{-2z}}\bigg]
\]
or, after having put
\[
e^{-z}=y,~~~ \text{whence},~~~ z=\log{\frac{1}{y}},~~~ e^{-z}\diff{d}z=-\diff{d}y
\]
\[
2. ~~~~ \int_0^{\infty} \frac{e^{au}-e^{-au}}{e^{\pi u}-e^{-\pi u}} \log{(x^2+u^2)} \diff{d}u= \int_0^1 \frac{\diff{d}y}{\log{\frac{1}{y}}}(\tan{\frac{1}{2}a}-\frac{2y^{x}\sin{a}}{1+2y\cos{a}+y^{2}})
\]
\[
(a<\pi)
\]
It can indeed easily be shown, if $a$ is in an arbitrary commensurable ratio to $\pi$, that is $a=\frac{m\pi}{n}$ ($m$ and $n$ integer numbers), that the integral on the right-hand side of equation (2.) can always be expressed in finite terms by the function $\Gamma$. For let, for the sake of brevity, be
\[
T= \int_0^1 \frac{\diff{d}y}{\log{\frac{1}{y}}}(\tan{\frac{1}{2}a}-\frac{2y^{x}\sin{a}}{1+2y\cos{a}+y^{2}}),
\]
after having differentiated with respect to $x$, it emerges
\[
3. ~~~~ \frac{\diff{d}T}{\diff{d}x}=2\int_0^1 \frac{2y^{x}\sin{a}\cdot\diff{d}y}{1+2y\cos{a}+y^{2}}
\]
But in the assumption $a=\frac{m\pi}{n}$ the integral is known \footnote{Legendre Exerc. book 2, \cite{8} pages 163-165 - We observe an error, which occurs at the cited place in Legendre's book, where he says about the formulas, mentioned by us in the text, that that one is valid, if $m$ is an odd number, this one on the other hand, if $m$ is an even number. But indeed it has to be prescribed, that that one or this one is valid, depending on, whether $m+n$ is odd or even. C.J. Malmstèn}
\[
\int_0^1 \frac{2y^{x}\sin{a}\cdot\diff{d}y}{1+2y\cos{a}+y^{2}}= \sum_{i=1}^{n-1}(-1)^{i-1}\sin{(ia)}[\frac{\diff{d}.\log{\Gamma{(\frac{x+n+i}{2n})}}}{\diff{d}x}-\frac{\diff{d}.\log{\Gamma{(\frac{x+i}{2n})}}}{\diff{d}x}]
\]
where $m+n$ is an odd number, and
\[
\int_0^1 \frac{2y^{x}\sin{a}\cdot\diff{d}y}{1+2y\cos{a}+y^{2}}= \sum_{i=1}^{\frac{1}{2}(n-1)}(-1)^{i-1}\sin{(ia)}[\frac{\diff{d}.\log{\Gamma{(\frac{x+n+i}{n})}}}{\diff{d}x}-\frac{\diff{d}.\log{\Gamma{(\frac{x+i}{n})}}}{\diff{d}x}]
\]
where $m+n$ is an even number; this, of course, multiplied by (3.), gives
\[
\diff{d}T=2\sum_{i=1}^{n-1}(-1)^{i-1}\sin{(ia)}[\frac{\diff{d}.\log{\Gamma{(\frac{x+n+i}{2n})}}}{\diff{d}x}-\frac{\diff{d}.\log{\Gamma{(\frac{x+i}{2n})}}}{\diff{d}x}]
\]
\[
(m+n= \text{odd number})
\]
\[
\diff{d}T=\sum_{i=1}^{\frac{1}{2}(n-1)}(-1)^{i-1}\sin{(ia)}[\frac{\diff{d}.\log{\Gamma{(\frac{x+n+i}{n})}}}{\diff{d}x}-\frac{\diff{d}.\log{\Gamma{(\frac{x+i}{n})}}}{\diff{d}x}]
\]
\[
(m+n= \text{even number})
\]
whence, after having integrated,\\

4.                 
$ \begin{cases}
T= C +2\sum_{i=1}^{n-1}(-1)^{i-1}\sin{(ia)}\log\bigg\{{\frac{\Gamma({\frac{x+n+i}{2n})}}{\Gamma({\frac{x+i}{2n})}}}\bigg\} \\[3mm]
~~~~~~~(m+n = \text{odd number})\\[3mm]
T= C^{\prime} +2\sum_{i=1}^{\frac{1}{2}(n-1)}(-1)^{i-1}\sin{(ia)}\log\bigg\{{\frac{\Gamma({\frac{x+n+i}{n})}}{\Gamma({\frac{x+i}{n})}}}\bigg\} \\[3mm]
~~~~~~~(m+n= \text{even number})
\end{cases}$\\[5mm]
Hence let us at first set $x=r$ and $x=s$ afterwards, by subtracting we will obtain this integral\\

5.
$ \begin{cases}
~~~~~~~~~~~\int_0^1 \frac{\diff{d}y}{\log{\frac{1}{y}}} \cdot \frac{y^r(1-y^{s-r})}{1+2y\cos{a}+y^2}\\[3mm]
=\csc{a}\sum_{i=1}^{n-1}(-1)^{i-1}\sin{(ia)}\log\bigg\{{\frac{\Gamma({\frac{s+n+i}{2n})}\cdot\Gamma({\frac{r+i}{2n})}}{\Gamma({\frac{r+n+i}{2n})}\cdot\Gamma({\frac{s+i}{2n})}}}\bigg\}\\[3mm]
~~~~~~~~~~(m+n= \text{odd number})\\[3mm]
~~~~~~~~~~~\int_0^1 \frac{\diff{d}y}{\log{\frac{1}{y}}} \cdot \frac{y^r(1-y^{s-r})}{1+2y\cos{a}+y^2}\\[3mm]
=\csc{a}\sum_{i=1}^{\frac{1}{2}(n-1)}(-1)^{i-1}\sin{(ia)}\log\bigg\{{\frac{\Gamma({\frac{s+n-i}{n})}\cdot\Gamma({\frac{r+i}{n})}}{\Gamma({\frac{r+n-i}{n})}\cdot\Gamma({\frac{s+i}{n})}}}\bigg\}\\[3mm]
~~~~~~~~~~(m+n= \text{even number})\\
\end{cases}$\\[5mm]

By the means of these formulas the values $C$ and $C^{\prime}$ can now indeed be determinend. For, if we set $r=0$ and $s=1$ and thereafter $r=1$, $s=2$, we, after having subtracted them from each other, find
\[
\int_0^1 \frac{\diff{d}y}{\log\frac{1}{y}}\cdot\frac{(1-2y+y^2)\sin{a}}{1+2y\cos{a}+y^2}= \sum_{i=1}^{n-1}(-1)^{i-1}\sin{ia}\log\bigg\{{\frac{(\Gamma{(\frac{n+i+1}{2n}}))^2\cdot\Gamma{(\frac{i+2}{2n}})\cdot\Gamma{(\frac{i}{2n}})}{(\Gamma({\frac{i+1}{2n})})^2\cdot\Gamma({\frac{n+i}{2n})}\cdot\Gamma({\frac{n+i+2}{2n})}}}\bigg\}
\]
\[
(m+n = \text{ odd number})
\]
\[
\int_0^1 \frac{\diff{d}y}{\log\frac{1}{y}}\cdot\frac{(1-2y+y^2)\sin{a}}{1+2y\cos{a}+y^2}= \sum_{i=1}^{n-1}(-1)^{i-1}\sin{ia}\log\bigg\{{\frac{(\Gamma{(\frac{n+1-i}{n}}))^2\cdot\Gamma{(\frac{i+2}{n}})\cdot\Gamma{(\frac{i}{n}})}{(\Gamma({\frac{i+1}{n})})^2\cdot\Gamma({\frac{n-i}{n})}\cdot\Gamma({\frac{n+2-i}{n})}}}\bigg\}
\]
\[
(m+n = \text{ even number})
\]

and from the formulas (4.), after having put $x=1$, it also emerges
\[
\int_0^1 \frac{\diff{d}y}{\log\frac{1}{y}}\cdot\frac{(1-2y+y^2)\sin{a}}{1+2y\cos{a}+y^2}
\]
\[
=(1+\cos{a})\bigg\{C+2\sum_{i=1}^{n-1}(-1)^{i-1}\sin{ia}\log\bigg\{{\frac{\Gamma{(\frac{n+i+1}{2n}})}{\Gamma{(\frac{i+1}{2n}})}}\bigg\}\bigg\}
\]
\[
(m+n= \text{odd number})
\]
\[
\int_0^1 \frac{\diff{d}y}{\log\frac{1}{y}}\cdot\frac{(1-2y+y^2)\sin{a}}{1+2y\cos{a}+y^2}
\]
\[
=(1+\cos{a})\bigg\{C^{\prime}+2\sum_{i=1}^{\frac{1}{2}(n-1)}(-1)^{i-1}\sin{ia}\log\bigg\{{\frac{\Gamma{(\frac{n+1-i}{n}})}{\Gamma{(\frac{i+1}{n}})}}\bigg\}\bigg\}
\]
\[
(m+n= \text{even number})
\]
which formulas, compared to each other, because we have in general
\[
2\sin{ma}\cos{a}=\sin{(m+1)a}+\sin{(m-1)a},
\]
give these values for $C$ and $C^{\prime}$ :

\begin{align*}
C&=\tan{\frac{1}{2}a}\cdot\log{2n}\\
C^{\prime}&=\tan{\frac{1}{2}a}\cdot\log{n}
\end{align*}
So, having substituted this in (4.), this becomes
\begin{align*}
T =& \tan{\frac{1}{2}a}\cdot\log{2n}+2\sum_{i=1}^{n-1}(-1)^{i-1}\sin{ia}\log\bigg\{{\frac{\Gamma{(\frac{x+n+i}{2n}})}{\Gamma{(\frac{x+i}{2n}})}}\bigg\}\\
&~~~~~~~(m+n=\text{ odd number})\\
T =& \tan{\frac{1}{2}a}\cdot\log{n}+2\sum_{i=1}^{\frac{1}{2}(n-1)}(-1)^{i-1}\sin{ia}\log\bigg\{{\frac{\Gamma{(\frac{x+n+i}{n}})}{\Gamma{(\frac{x+i}{n}})}}\bigg\}\\
&~~~~~~~(m+n=\text{ even number}),\\
\end{align*}
and if we, for the sake of brevity, set 
\[
6. ~~~~ L(a,x)=\int_0^{\infty} \frac{e^{au}-e^{-au}}{e^{\pi u}-e^{-\pi u}} \log{(x^2+u^2)} \diff{d}u
\]
if follows from formula (2.), after having put $a=\frac{m\pi}{n}$, where $m<n$,\\

7. 
$\begin{cases}
L(a,x)= \tan{\frac{1}{2}a}\cdot\log{2n}+2\sum_{i=1}^{n-1}(-1)^{i-1}\sin{ia}\log\bigg\{{\frac{\Gamma{(\frac{x+n+i}{2n}})}{\Gamma{(\frac{x+i}{2n}})}}\bigg\}\\[3mm]
~~~~~~~~~~~~~~~~~~~~~~(m+n=\text{ odd number})\\[3mm]
L(a,x)=\tan{\frac{1}{2}a}\cdot\log{n}+2\sum_{i=1}^{\frac{1}{2}(n-1)}(-1)^{i-1}\sin{ia}\log\bigg\{{\frac{\Gamma{(\frac{x+n-i}{n}})}{\Gamma{(\frac{x+i}{n}})}}\bigg\}\\[3mm]
~~~~~~~~~~~~~~~~~~~~~~(m+n=\text{ even number})
\end{cases} $\\[10mm]

For $x=1$ and $a=\frac{1}{2}\pi$, whence $m=1$, $n=2$, you get from the first
\[
8.~~~~~~ \int_0^{\infty} \frac{\log{(1+u^2)}}{e^{\frac{1}{2}\pi u}-e^{-\frac{1}{2}\pi u}}\diff{d}u=\log{(\frac{4}{\pi})}
\]
\section*{~~~~~~~~~~~~~~~~~~~~~~~~~~~~~~~~~§2}
But to deduce certain certain new definite integrals from the formulas (7.), it is at first neccessary to demonstrate, that these formulas are also valid for $x=0$; That this is indeed the case, if suffices to constitute, that formula (1.), whence those have their origin, also remains correct for $x=0$.\\
The following formulas are known

\begin{align*}
\int_0^{\omega} e^{-xz}\cos{uz}\diff{d}z&=\frac{e^{-\omega z}(u\sin{u\omega}-x\cos{u \omega})}{x^2+u^2}+\frac{x}{x^2+u^2}\\
\int_0^{\omega} e^{-uz}\sin{xz}\diff{d}z&=-\frac{e^{-\omega u}(u\sin{\omega x}-x\cos{ \omega x})}{x^2+u^2}+\frac{x}{x^2+u^2}
\end{align*}
whence, after having subtracted, we will have
\[
\int_0^{\omega} \diff{d}z(e^{-xz}\cos{uz}-e^{-uz}\sin{xz})
\]
\[
=\frac{e^{-\omega z}(u\sin{u\omega}-x\cos{u \omega})}{x^2+u^2}+\frac{e^{-\omega u}(u\sin{\omega x}-x\cos{ \omega x})}{x^2+u^2}
\]
Now, if, after having multiplied by $\diff{d}x$ on both sides, we take the integral from $x=0$ to $x=u$, it is
\[
\int_0^{\omega} \frac{\diff{d}z}{z}(\cos{uz}-e^{-uz})
\]
\[
=e^{-\omega u} \int_0^u\frac{(u\sin{\omega x}+x\cos{\omega x})}{x^2+u^2}\diff{d}x ~~\text{and} ~~ \int_0^u \frac{e^{-\omega x}(u\sin{u \omega}-x\cos{u\omega})}{x^2+u^2}\diff{d}x
\]
It indeed easily becomes clear, that for the values $\omega$, increasing indefinitely, the integrals
\[
e^{-\omega u} \int_0^u\frac{(u\sin{\omega x}+x\cos{\omega x})}{x^2+u^2}\diff{d}x ~~\text{and} ~~ \int_0^u \frac{e^{-\omega x}(u\sin{u \omega}-x\cos{u\omega})}{x^2+u^2}\diff{d}x
\]
converge to zero indefinitely, whence it follows, that
\[
\lim \int_0^{\omega} \frac{\diff{d}z}{z}(\cos{uz}-e^{-uz})=0 ~~~~ (\omega = \infty)
\]
that is
\[
\int_0^{\infty} \frac{\diff{d}z}{z}(\cos{uz}-e^{-uz})=0 ~~~~ (\omega = \infty)
\]
one will obtain, if it is subtracted from that known one
\[
\int_0^{\infty} \frac{e^{-z}-e^{uz}}{z}\diff{d}z=\log{u},
\]
this formula 
\[
9.~~~~~~~ \int_0^{\infty} \frac{e^{-z}-\cos{uz}}{z}\diff{d}z=\log{u},
\]
whence it becomes obvious, that formula (1.) along with all derived from it are valid for $x=0$.\\

After having demonstrated these things this way, now let be $x=0$ in (7.), and let us put
\[
e^{\frac{\pi u}{n}}=y, ~~\text{whence}~~ u=\frac{n}{\pi}\log{y}
\]
because $a=\frac{m \pi}{n}$, you will obtain
\[
\int_0^{\infty} \frac{y^{m-1}-y^{-m-1}}{y^n-y^{-n}}\{\log{(\frac{n}{\pi})}+\log{(\log{y})}\}
\]
\[
=\frac{\pi}{2n}\cdot\tan{\frac{m \pi}{2n}}\cdot\log{2n}+\frac{\pi}{n}\cdot\sum_{i=1}^{n-1}(-1)^{i-1}\sin{\frac{im\pi}{n}}\log\bigg\{{\frac{\Gamma{(\frac{n+i}{2n}})}{\Gamma{(\frac{i}{2n}})}}\bigg\}
\]
\[
(m+n= \text{odd number})
\]
\[
\int_0^{\infty} \frac{y^{m-1}-y^{-m-1}}{y^n-y^{-n}}\{\log{(\frac{n}{\pi})}+\log{(\log{y})}\}
\]
\[
=\frac{\pi}{2n}\cdot\tan{\frac{m \pi}{2n}}\cdot\log{n}+\frac{\pi}{n}\cdot\sum_{i=1}^{\frac{1}{2}(n-1)}(-1)^{i-1}\sin{\frac{im\pi}{n}}\log\bigg\{{\frac{\Gamma{(1-\frac{1}{n}})}{\Gamma{(\frac{i}{n}})}}\bigg\}
\]
\[
(m+n= \text{even number})
\]
and because this formula holds 
\[
\int_0^{\infty} \frac{y^{m-1}-y^{-m-1}}{y^n-y^{-n}}\diff{d}y= -\int_0^{1} \frac{z^{\frac{m}{2n}-\frac{1}{2}}-z^{-\frac{m}{2n}-\frac{1}{2}}}{1-z}\diff{d}z=\frac{\pi}{2n}\tan{\frac{m \pi}{2n}}
\]
also\\[2mm]

10. $\begin{cases}
~~~~~~~~~~~~~~~~~\int_0^{\infty} \frac{y^{m-1}-y^{-m-1}}{y^n-y^{-n}}\log{(\log{y})}\diff{d}y\\[5mm]
=\frac{\pi}{2n}\cdot\tan{\frac{m \pi}{2n}}\cdot\log{2\pi}+\frac{\pi}{n}\cdot\sum_{i=1}^{n-1}(-1)^{i-1}\sin{\frac{im\pi}{n}}\log\bigg\{{\frac{\Gamma{(\frac{n+i}{2n}})}{\Gamma{(\frac{i}{2n}})}}\bigg\}\\[5mm]
~~~~~~~~~~~~~~~~~~~~~~~~(m+n= \text{odd number})\\[5mm]
~~~~~~~~~~~~~~~~~\int_0^{\infty} \frac{y^{m-1}-y^{-m-1}}{y^n-y^{-n}}\log{(\log{y})}\diff{d}y\\[5mm]
=\frac{\pi}{2n}\cdot\tan{\frac{m \pi}{2n}}\cdot\log{\pi}+\frac{\pi}{n}\cdot\sum_{i=1}^{\frac{1}{2}(n-1)}(-1)^{i-1}\sin{\frac{im\pi}{n}}\log\bigg\{{\frac{\Gamma{(1-\frac{1}{n}})}{\Gamma{(\frac{i}{n}})}}\bigg\}\\[5mm]
~~~~~~~~~~~~~~~~~~~~~~~~(m+n= \text{even number})
\end{cases}$\\[5mm]

Hence it follows for m=1\\[2mm]

11. $\begin{cases}
~~~~~~~~~~~~~~~~~\int_0^{\infty} \frac{y^{n-2}\log{(\log{y})}}{1+y^2+y^4+\cdots+y^{2(n-1)}}\diff{d}y\\[5mm]
=\frac{\pi}{2n}\cdot\tan{\frac{m \pi}{2n}}\cdot\log{2\pi}+\frac{\pi}{n}\cdot\sum_{i=1}^{n-1}(-1)^{i-1}\sin{\frac{i\pi}{n}}\log\bigg\{{\frac{\Gamma{(\frac{n+i}{2n}})}{\Gamma{(\frac{i}{2n}})}}\bigg\}\\[5mm]
~~~~~~~~~~~~~~~~~~~~~~~~(n= \text{even number})\\[5mm]
~~~~~~~~~~~~~~~~~\int_0^{\infty} \frac{y^{n-2}\log{(\log{y})}}{1+y^2+y^4+\cdots+y^{2(n-1)}}\diff{d}y\\[5mm]
=\frac{\pi}{2n}\cdot\tan{\frac{m \pi}{2n}}\cdot\log{\pi}+\frac{\pi}{n}\cdot\sum_{i=1}^{\frac{1}{2}(n-1)}(-1)^{i-1}\sin{\frac{i\pi}{n}}\log\bigg\{{\frac{\Gamma{(1-\frac{1}{n}})}{\Gamma{(\frac{i}{n}})}}\bigg\}\\[5mm]
~~~~~~~~~~~~~~~~~~~~~~~~(n= \text{odd number})
\end{cases}$\\[5mm]

For the sake of an example, after having put $n=2$ in that one and $n=3$ in this one, we will have, if in this one $y^2$ is changed into $y$,\\

12\footnote{The first of the following integrals is nowadays called Vardi's integral, after Vardi, who considered it in \cite{10}, but its origin is this paper, written by Malmstén.}  $~~~~~~~~\begin {cases}
\int_0^{\infty} \frac{\log{(\log{y})}\diff{d}y}{1+y^2}=\frac{1}{2}\pi\log\bigg\{{\frac{(2\pi)^{\frac{1}{2}}\Gamma{(\frac{3}{4}})}{\Gamma{(\frac{1}{4}})}}\bigg\}\\[5mm]
\int_0^{\infty} \frac{\log{(\log{y})}\diff{d}y}{1+y+y^2}=\frac{\pi}{\sqrt{3}}\log\bigg\{{\frac{(2\pi)^{\frac{1}{3}}\Gamma{(\frac{2}{3}})}{\Gamma{(\frac{1}{3}})}}\bigg\}\\
\end{cases} $
\section*{~~~~~~~~~~~~~~~~~~~~~~~~~~~~~~~~~§3}

From the first §. it became clear already, that that transcendaental, we denoted by $L(a,x)$, can always be expressed in finite terms by the means of $\Gamma$, as long as $a$ is in a commensurable ratio to $\pi$. Now in this §. we want to deal with certain of its remarkable properties. So let be $a=\frac{1}{2}\pi$ in the first of the formulas (7.), then it will be
\[
13. ~~~ L(\frac{1}{2}\pi,x)= 2\log\bigg\{{\frac{2\Gamma({\frac{x+3}{4}})}{\Gamma({\frac{x+1}{4}})}}\bigg\}
\]
whence, if one put $x+2$ in the place of $x$, we find by adding
\[
14. ~~ L(\frac{1}{2}\pi,x+2)+L(\frac{1}{2}\pi,x)=2\log{(x+1)}
\]
But if we substitute $2-x$ in the place of $x$ in (13.), from the known relation\footnote{Euler proved his famous reflection formula in \cite{3}.}
\[
\Gamma{(a)} \Gamma{(1-a)}=\frac{\pi}{\sin{a \pi}}
\]  
it will be in the same way by adding
\[
15. ~~ L(\frac{1}{2}\pi,x)+L(\frac{1}{2}\pi,2-x)=2\log{[(x+1)\cot{\frac{(x+1)\pi}{4})]}}
\]
What, for $x=1$, renders formula (8.)
\[
16. ~~ L(\frac{1}{2}\pi,1)=\log{(\frac{4}{\pi})}
\]
From the formulas (14. and 15.) we know the function $L(\frac{1}{2}\pi,x)$ for any arbitrary value of $x$, if it is only known throughout the period from $x=0$ to $x=1$. Moreover the formulas (14. and 16.) teach, that $L(\frac{1}{2}\pi,x)$, for $x=$ any arbitrary integer, can be expressed finitely by means of logarithms and $\pi$.\\

Now let be $a=\frac{2}{3}\pi$ in the first of the formulas (7.), then it will be
\[
L(\frac{2}{3}\pi,x)=2 \sin{\frac{1}{3}\pi}\log\bigg\{{\frac{6\cdot\Gamma({\frac{x+4}{6}})\Gamma({\frac{x+5}{6}})}{\Gamma({\frac{x+1}{6}})\Gamma({\frac{x+2}{6}})}}\bigg\}
\]
whence in almost the same way as above we will obtain\\[2mm]

17. $\begin{cases}
L(\frac{2}{3}\pi,x+3)+L(\frac{2}{3}\pi,x)= 2\sin{\frac{1}{3}\pi}\log{[(x+2)(x+1)]}\\[5mm]
~~~~~~~~~~L(\frac{2}{3}\pi,x)+L(\frac{2}{3}\pi,3-x)\\[5mm]
=2\sin{\frac{1}{3}\pi}\log{[(x-2)(x-1)\tan{\frac{(x+2)\pi}{6}}\tan{\frac{(x+1)\pi}{6}}]}\\
\end{cases}$\\[5mm]

and from the secon for $x=\frac{3}{2}$
\[
18.~~~~L(\frac{2}{3}\pi,\frac{3}{2})=2\sin{\frac{1}{3}\pi}\log({\frac{1}{2}\cot{\frac{\pi}{12}})}
\]

and from the formulas (17.) it follows, that the function $L(\frac{2}{3}\pi,x)$ is known for all values of $x$, if we only know it for any arbitrary value between $x=0$ and $x=\frac{3}{2}$, moreover the first of these formulas along with (18.) teaches, that $L(\frac{2}{3}\pi,x)$ can be expressed in finite terms by means of logarithms and trigonometrical functions for $x=\frac{1}{2}((2i+1)\cdot3)$, while $i$ denotes any arbitrary integer number.\\

But we can indeed find a far more general relation, from which the preceeding can be derived as special cases. For, where $a=\frac{m\pi}{n}$ and $m+n=$ odd number, we also have
\[
19.~~~ \tan{\frac{1}{2}a}= \sum_{i=1}^{n-1} (-1)^{i-1}\sin{ia},
\]
from where the first of the formulas (7.) can be transformed into this form:
\[
20.~~~L(a,x) =2\sum_{i=1}^{n-1} (-1)^{i-1}\sin{ia}\log\bigg\{{\frac{(2n)^{\frac{1}{2}}\cdot\Gamma({\frac{x+n+i}{2n}})}{\Gamma({\frac{x+i}{2n}})}}\bigg\}.
\]
\[
(m+n = \text{odd number})
\]
Let us substitute $x+n$ for $x$ here; then we will have by adding
\[
21.~~~ L(a,x+n)+L(a,x)= 2\sum_{i=1}^{n-1} (-1)^{i-1}\sin{ia}\log{(x+i)}
\]
\[
(m+n = \text{odd number})
\]
further, if in (20.) $n-x$ is substituted in the place $x$, it also emerges, by adding,
\[
L(a,x)+L(a,n-x)=2\sum_{i=1}^{n-1} (-1)^{i-1}\sin{ia}\log\bigg\{{\frac{(2n)^{\frac{1}{2}}\cdot\Gamma({\frac{1}{2}+\frac{x+i}{2n}})}{\Gamma({\frac{x+i}{2n}})}}\bigg\}
\]
\[
+2\sum_{i=1}^{n-1} (-1)^{i-1}\sin{ia}\log\bigg\{{\frac{(2n)^{\frac{1}{2}}\cdot\Gamma({1+\frac{i-x}{2n}})}{\Gamma({\frac{1}{2}+\frac{i-x}{2n}})}}\bigg\}
\]
or, because we have
\[
~\sum_{i=1}^{n-1} (-1)^{i-1}\sin{ia}\log\bigg\{{\frac{(2n)^{\frac{1}{2}}\cdot\Gamma({1+\frac{i-x}{2n}})}{\Gamma({\frac{1}{2}+\frac{i-x}{2n}})}}\bigg\}
\]
\[
=\sum_{i=1}^{n-1} (-1)^{i-1}\sin{ia}\log\bigg\{{\frac{(2n)^{\frac{1}{2}}\cdot\Gamma({\frac{3}{2}-\frac{x+i}{2n}})}{\Gamma({1-\frac{x+i}{2n}})}}\bigg\},
\]
also 
\[
L(a,x)+L(a,n-x)
\]
\[
=2\sum_{i=1}^{n-1} (-1)^{i-1}\sin{ia}\log\bigg\{{\frac{2n\cdot\Gamma({\frac{3}{2}-\frac{x+i}{2n}})\cdot\Gamma({\frac{1}{2}+\frac{x+i}{2n}})}{\Gamma({\frac{x+i}{2n}})\cdot\Gamma({1-\frac{x+i}{2n}})}}\bigg\}.
\]
But it is
\[
\Gamma({\frac{3}{2}-\frac{x+i}{2n}})\cdot\Gamma({\frac{1}{2}+\frac{x+i}{2n}})=\frac{n-x-i}{2n}\cdot\frac{\pi}{\sin{\frac{(x+1)\pi}{2n}}},
\]
\[
\Gamma({\frac{x+i}{2n}})\cdot\Gamma({1-\frac{x+i}{2n}})=\frac{\pi}{\sin{\frac{(x+1)\pi}{2n}}},
\]
whence it follows
\[
L(a,x)+L(a,n-x)=2\sum_{i=1}^{n-1} (-1)^{i-1}\sin{ia}\log{[(n-x-i)\tan{\frac{1}{2n}(x+i)\pi}]},
\]
or at least, after having put $n-x$ in the place of $x$
\[
22. ~~~L(a,x)+L(a,n-x)=2\sum_{i=1}^{n-1} (-1)^{i-1}\sin{ia}\log{[(x-i)\cot{\frac{(x-i)\pi}{2n}}}]
\]
\[
(m+n= \text{odd number})
\]
So it is clear from the formulas (21. and 22.), that the function $L(a,x)$ is known for each value of $x$, if it is only known throughout the whole period from $x=0$ to $x=\frac{1}{2}n$, while $m+n$ is an odd number.\\

After having put $x=\frac{1}{2}n$ in (22.), we have the formula
\[
23.~~~ L(a, \frac{1}{2}n)=\sum_{i=1}^{n-1} (-1)^{i-1}\sin{ia}\log{[(\frac{1}{2}n-i)\cot({\frac{\pi}{4}-\frac{i \pi}{2n}})}],
\]
which, of course, for $i=\frac{1}{2}n$ presents the expression $\log({0\cdot \infty})$; nevertheless its true value is easily found to be $\log{(\frac{2n}{\pi})}$. This formula along with (21.) teaches, that the function $L(a,x)$, if $m+n$ is an odd number, can be expressed in finite terms by means of logarithms and trigonometrical functions, after having put $x=\frac{1}{2}n(2i+1)$, where $i$ denotes a certain integer number.\\

Almost likewise we can derive analogous relations from the second of the formula (7.). For substitute $x+n$ in the place $x$, then it will be
\[
L(a,x+n)= \tan{\frac{1}{2}a}\cdot\log{n}+2\sum_{i=1}^{\frac{1}{2}(n-1)}(-1)^{i-1}\sin{ia}\log\bigg\{{\frac{\Gamma{(\frac{2n+x-i}{n}})}{\Gamma{(\frac{n+x+i}{n}})}}\bigg\},
\]
whence by subtracting
\[
24. ~~~ L(a,x+n)-L(a,x)=2\sum_{i=1}^{\frac{1}{2}(n-1)}(-1)^{i-1}\sin{ia}\log\bigg\{{\frac{x+n-i}{x+i}}\bigg\}
\]
\[
(m+n = \text{even number})
\]
Afterall, if in the same formula (7.) one puts $n-x$ for $x$, it becomes
\[
L(a,n-x)= \tan{\frac{1}{2}a}\cdot\log{n}+2\sum_{i=1}^{\frac{1}{2}(n-1)}(-1)^{i-1}\sin{ia}\log\bigg\{{\frac{\Gamma{(\frac{2n-x-i}{n}})}{\Gamma{(\frac{n+i-x}{n}})}}\bigg\},
\]
whence by subtraction from a known relation of the function $\Gamma$
\[
25. ~~~ L(a,x)-L(a,n-x)=2\sum_{i=1}^{\frac{1}{2}(n-1)}(-1)^{i-1}\sin{ia}\log\bigg\{{\frac{x+n-i}{x+i}\cdot\frac{\sin{\frac{(x+i)\pi}{n}}}{\sin{\frac{(x-i)\pi}{n}}}}\bigg\}
\]
\[
(m+n = \text{even number})
\]
By means of the formula (24. and 25.) it follows, that, whilw $m+n$ is an even number, that we know the function $L(a,x)$ all over for each value $x$, if it is only known throughout the period from $x=0$ to $x=\frac{1}{2}$.\\

Now let us successively put in the place of $x$:
\[
x,~~~x+\frac{2n}{r}, ~~~x+\frac{4n}{r}, ~~~x+\frac{6n}{r},\cdots, ~~~x+\frac{(r-1)\cdot2n}{r}
\]
in the first of the formulas (7.), and 
\[
x,~~~x+\frac{n}{r}, ~~~x+\frac{2n}{r}, ~~~x+\frac{3n}{r},\cdots, ~~~x+\frac{(r-1)\cdot n}{r}
\]
in the second; then by means of the known\footnote{This is the Gau\ss ian multiplication formula. Gau\ss{}  proved it in \cite{6}, but Malmstèn most likely got it from Legendre's book \cite{8}.} formula
\[
\Gamma{(y)}\cdot\Gamma{(y+\frac{1}{r})}\cdot\Gamma{(y+\frac{2}{r})}\cdots\cdot\Gamma{(y+\frac{r-1}{r})}= \Gamma{(ry)}\cdot(2 \pi)^{\frac{1}{2}(r-1)}\cdot r^{\frac{1}{2}-ry}
\]
we find by  summing\\[2mm]
26.
$\begin{cases}
~~~L(a,x)+L(a, x+\frac{2n}{r})+L(a, x+\frac{4n}{r})+\cdots+L(a, x+\frac{(r-1)\cdot 2n}{r})\\[5mm]
=r\tan{\frac{1}{2}a}\cdot\log{2n}+2\sum_{i=1}^{n-1}(-1)^{i-1}\sin{ia}\log\bigg\{{\frac{\Gamma({\frac{1}{2}r+\frac{r(x+i}{n}})}{r^{\frac{1}{2}r}\Gamma{(\frac{r(x+i)}{2n}})}}\bigg\}\\[5mm]
~~~~~~~~~~~~~~~~~~(m+n= \text{odd number})\\[5mm]
~~~L(a,x)+L(a, x+\frac{n}{r})+L(a, x+\frac{2n}{r})+\cdots+L(a, x+\frac{(r-1)\cdot n}{r})\\[5mm]
=r\tan{\frac{1}{2}a}\cdot\log{n}+2\sum_{i=1}^{\frac{1}{2}(n-1)}(-1)^{i-1}\sin{ia}\log\bigg\{{\frac{\Gamma({\frac{1}{2}+\frac{r(x-i}{n}})}{r^{r-\frac{2ri}{n}}\Gamma{(\frac{r(x+i)}{n}})}}\bigg\}\\[5mm]
~~~~~~~~~~~~~~~~~~(m+n= \text{even number})\\
\end{cases}$\\[2mm]

It easily becomes clear, that the first of these sums can be expressed by the means of logarithms, while $r$ is an even number.\\

\section*{~~~~~~~~~~~~~~~~~~~~~~~~~~~~~~~~~§4}
This formula is known\footnote{This formula is due to Euler \cite{4}}
\[
27.~~~ \int_0^{\infty}z^{s-1}e^{-xz}\cos{(uz)}\diff{d}z= \frac{\Gamma{(s)}}{(x^2+u^2)^{\frac{1}{2}s}}\cdot\cos{(s\arctan{\frac{u}{x}})}
\]
which, while $1>s>0$  also legitimate for $x=0$, gives
\[
28.~~~
\int_0^{\infty}z^{s-1}\cos{(uz)}\diff{d}z= \frac{\cos{(\frac{s\pi}{2})}\cdot\Gamma{(s)}}{u^{s}}.
\]
Let us multyply (27.) on both sides by
\[
\frac{e^{au}-e^{-au}}{e^{\pi u}-e^{-\pi u}}\cdot\diff{d}u ~~~~[a<\pi];
\]
after having integrated between $u=0$ and $u=\infty$, it becomes
\[
\int_0^{\infty} \frac{e^{au}-e^{-au}}{e^{\pi u}-e^{-\pi u}}\cdot \frac{\cos{(s\arctan{\frac{u}{x}})}}{(x^2+u^2)^{\frac{s}{2}}}\cdot \diff{d}u
\]
\[
=\frac{1}{\Gamma{(s)}} \cdot \int_0^{\infty} e^{-xz}\cdot z^{s-1}\diff{d}z\cdot \int_0^{\infty} \frac{e^{au}-e^{-au}}{e^{\pi u}-e^{-\pi u}} \cdot \cos{uz} \cdot \diff{d}u
\]
\[
\frac{\sin{a}}{\Gamma{(s)}}\cdot \int_0^{\infty} \frac{e^{-xz}z^{s-1} \cdot e^{-z} \diff{d}z}{1+2e^{-z}\cos{a}+e^{-2z}},
\]
whence, by putting $e^{-z}=y$ on the right-hand side, it emerges
\[
29.~~~ \int_0^{\infty} \frac{e^{au}-e^{-au}}{e^{\pi u}-e^{-\pi u}}\cdot \frac{\cos{(s\arctan{\frac{u}{x}})}}{(x^2+u^2)^{\frac{s}{2}}}\cdot \diff{d}u = \frac{\sin{a}}{\Gamma{(s)}}\cdot \int_0^{1} \frac{y^x(\log\frac{1}{y})^{s-1}\diff{d}y}{1+2y\cos{a}+y^2}.
\]
and for $x=0$, if only $1>s>0$,
\[
30.~~~ \int_0^{\infty} \frac{e^{au}-e^{-au}}{e^{\pi u}-e^{-\pi u}}\cdot \frac{\diff{d}u}{u^s}= \frac{\sin{a}}{\cos{\frac{1}{2}s \pi}\cdot \Gamma{(s)}}\cdot \int_0^{1} \frac{(\log\frac{1}{y})^{s-1}\diff{d}y}{1+2y\cos{a}+y^2}.
\]
From this formula, if one puts $a=\frac{m \pi}{n}$ and $e^{-\frac{\pi u}{n}}=y$, whence it is
\[
u=\frac{n}{\pi}\cdot \log{\frac{1}{y}}~~ \text{and}~~ \diff{d}u= -\frac{n}{\pi}\cdot \frac{\diff{d}y}{y},
\]
one finds, after having done the transformation, 
\[
31.~~ \int_0^1 \frac{y^{m-1}+y^{-m-1}}{y^n-y^{-n}}\cdot \frac{\diff{d}y}{(\log{\frac{1}{y}})^s}= \frac{(\frac{\pi}{n})^{1-s}\cdot \sin{\frac{m \pi}{n}}}{\cos{\frac{1}{2}s \pi}\cdot\Gamma{(s)}}\cdot \int_0^{1} \frac{(\log\frac{1}{y})^{s-1}\diff{d}y}{1+2y\cos{\frac{m \pi}{n}}+y^2},
\]
and for $m=1$,
\[
32.~~ \int_0^1 \frac{y^{n-2}\cdot \diff{d}y}{1+y^2+\cdots+y^{2(n-1)}}\cdot \frac{1}{(\log{\frac{1}{y}})^s}= \frac{(\frac{\pi}{n})^{1-s}\cdot \sin{\frac{ \pi}{n}}}{\cos{\frac{1}{2}s \pi}\cdot\Gamma{(s)}}\cdot \int_0^{1} \frac{(\log\frac{1}{y})^{s-1}\diff{d}y}{1+2y\cos{\frac{ \pi}{n}}+y^2},
\]
Now let us call
\[
33.~~ G(s)= \int_0^1 \frac{(\log\frac{1}{y})^{s-1}\diff{d}y}{1+y+y^2}~~,~~ G_1(s)= \int_0^1 \frac{(\log\frac{1}{y})^{s-1}\diff{d}y}{1+y^2};
\]
formula (32.) will give, after having put $n=3$ and having changed $y^2$ into $y$:
\[
34.~~ G(1-s)=\frac{(\frac{2}{3}\pi)^{1-s}\cdot \sin{\frac{1}{3}\pi}}{\cos{\frac{1}{2}s \pi}\cdot\Gamma{(s)}}\cdot G(s),
\]
and one concludes the same for $n=2$ immediately
\[
35.~~ G_1(1-s)=\frac{(\frac{1}{2}\pi)^{1-s}}{\cos{\frac{1}{2}s \pi}\cdot\Gamma{(s)}}\cdot G_1(s),
\]
Glow and behold the simple and remarkable equations, which connect the functions $G(s)$ and $G_1(s)$ and their complementary ones $G(1-s)$ and $G_1(1-s)$ to each other; in this respect they are analogous to this known formula of the function $\Gamma$
\[
\Gamma{(a)}\cdot \Gamma{(1-a)} = \frac{\pi}{\sin{a \pi}}
\]
From the formula (34. and 35) we will obtain by taking logarithms 
\begin{align*}
&\log{G(1-s)}+\log{G(s)}~~=(1-s)\log{\frac{2}{3}\pi}+\log{\sin{\frac{1}{3}\pi}}-\log{\cos{\frac{1}{2}s \pi}}-\log{\Gamma{(s)}}\\
&\log{G_1(1-s)}+\log{G_1(s)}=(1-s)\log{\frac{1}{2}\pi}-\log{\cos{\frac{1}{2}s \pi}}-\log{\Gamma{(s)}}
\end{align*}
whence, if we put for the sake of brevity\\[2mm]
36.
$ \begin{cases}
F(s)~=\frac{\diff{d}. G(s)}{\diff{d}s}=\int_0^1 \frac{(\log\frac{1}{y})^{s-1}\log{(\log{\frac{1}{y}})}\diff{d}y}{1+y+y^2}\\[5mm]
F_1(s)=\frac{\diff{d}. G_1(s)}{\diff{d}s}=\int_0^1 \frac{(\log\frac{1}{y})^{s-1}\log{(\log{\frac{1}{y}})}\diff{d}y}{1+y^2}
\end{cases}$\\[5mm]

by differentiating we will have these new relations\\[2mm]
37.
$ \begin{cases}
\frac{F(s)}{G(s)}+\frac{F(1-s)}{G(1-s)}~~=\log{\frac{2}{3}\pi}+Z^{\prime}(s)-\frac{1}{2}\pi \tan{\frac{1}{2}s \pi},\\[5mm]
\frac{F_1(s)}{G_1(s)}+\frac{F_1(1-s)}{G_1(1-s)}=\log{\frac{1}{2}\pi}+Z^{\prime}(s)-\frac{1}{2}\pi \tan{\frac{1}{2}s \pi},
\end{cases}$\\[5mm]

if, along with Legendre, we denote $\frac{\diff{d}. \log{\Gamma{(s)}}}{\diff{d}s}$ by $Z^{\prime}(s)$.\\

Let us suppose $s=\frac{1}{2}$ in (37.); then it will be
\begin{align*}
&\int_0^1 \frac{\log{(\log{\frac{1}{y})}}}{1+y+y^2}\cdot \frac{\diff{d}y}{\sqrt{\log{\frac{1}{y}}}}=\frac{1}{2}(\log{\frac{\pi}{6}}-\frac{1}{2}\pi-C)\cdot \int_0^1 \frac{1}{1+y+y^2}\cdot \frac{\diff{d}y}{\sqrt{\log{\frac{1}{y}}}},\\
&\int_0^1 \frac{\log{(\log{\frac{1}{y})}}}{1+y^2}\cdot \frac{\diff{d}y}{\sqrt{\log{\frac{1}{y}}}}=\frac{1}{2}(\log{\frac{\pi}{8}}-\frac{1}{2}\pi-C)\cdot \int_0^1 \frac{1}{1+y^2}\cdot \frac{\diff{d}y}{\sqrt{\log{\frac{1}{y}}}},
\end{align*}

or by putting $\log{\frac{1}{y}}=x$:\\[5mm]
38. $\begin{cases}
\int_0^{\infty} \frac{\log{x}}{e^x+1+e^{-x}} \cdot \frac{\diff{d}x}{\sqrt{x}}=(\log{\frac{\pi}{6}}-\frac{1}{2}\pi-C)\cdot \int_0^{\infty} \frac{\diff{d}x}{e^x+1+e^{-x}} \cdot \frac{1}{\sqrt{x}},\\[5mm]
\int_0^{\infty} \frac{\log{x}}{e^x+e^{-x}} \cdot \frac{\diff{d}x}{\sqrt{x}}~~~=(\log{\frac{\pi}{8}}-\frac{1}{2}\pi-C)\cdot \int_0^{\infty} \frac{\diff{d}x}{e^x+e^{-x}} \cdot \frac{1}{\sqrt{x}},
\end{cases}$\\[2mm]

while $C=-Z^{\prime}(1)$ is the known Eulerian constant $0,577216 \cdots$.These formulas, which we have not seen proposed up to now, do not seem unworthy of the Geometers' attention.\\

From the formulas (38.) we will now deduce two relations, which in the transformations of series seem to be worth of the highest places. Of course, while it is in identical manner\\[5mm]
39. $\begin{cases}
\frac{\sin{\frac{1}{3}\pi}}{e^x+1-e^{-x}}= \sum_{i=1}^n (-1)^{i-1}\cdot e^{ix} \cdot \sin{\frac{1}{3}i \pi} +\frac{(-1)^n \cdot e^{-nx}(\sin{\frac{1}{3}(n+1) \pi}+e^{-x}\sin{\frac{1}{3}n \pi})}{1+e^{-x}+e^{-2x}}\\[5mm]
\frac{1}{e^x-e^{-x}}~~~= \sum_{i=0}^{n-1} (-1)^{i}\cdot e^{-(2i+1)x} +\frac{(-1)^n  e^{-2nx}}{e^{x}+e^{-x}}
\end{cases}$\\[5mm]

from the known formula
\[
40.~~ \int_0^{\infty} e^{-kx} \cdot \log{x} \cdot \frac{\diff{d}x}{\sqrt{x}}= -\frac{\sqrt{\pi}}{\sqrt{k}}(\log{k}+2\log{2}+C)
\]
we will have
\[
\sin{\frac{1}{3}\pi} \cdot \int_0^{\infty} \frac{\log{x}}{e^x+1+e^{-x}} \cdot \frac{\diff{d}x}{\sqrt{x}}
\]
\[
=\sqrt{\pi} \cdot \sum_{i=1}^n (-1)^i \cdot \sin{\frac{1}{3} i \pi} \cdot \frac{\log{i}+2\log{2}+C}{\sqrt{i}}
\]
\[
+(-1)^n\{P(n) \cdot \sin{\frac{1}{3}(n+1)\pi} +P(n+1)\cdot \sin{\frac{1}{3}(n+1)\pi} \},
\]
\[
\int_0^{\infty} \frac{\log{x}}{e^x+e^{-x}} \cdot \frac{\diff{d}x}{\sqrt{x}}=\sqrt{\pi} \cdot \sum_{i=0}^{n-1} (-1)^{i+1}  \cdot \frac{\log{(2i+1)}+2\log{2}+C}{\sqrt{2i+1}}+(-1)^n \cdot Q(n),
\]
where we put, for the sake of brevity,
\[
P(n)=\int_0^{\infty} \frac{e^{-nx}\log{x}}{e^x+1+e^{-x}} \cdot \frac{\diff{d}x}{\sqrt{x}}=\theta \cdot \int_0^{\infty} \frac{e^{-nx} \log{x} \cdot \diff{d}x}{\sqrt{x}}= -\frac{\theta \cdot \sqrt{\pi}}{\sqrt{n}}(\log{n}+2 \log{2}+C),
\]
\[
Q(n)=\int_0^{\infty} \frac{e^{-2nx}\log{x}}{e^x+e^{-x}} \cdot \frac{\diff{d}x}{\sqrt{x}}=\theta_1 \cdot \int_0^{\infty} \frac{e^{-2nx} \log{x} \cdot \diff{d}x}{\sqrt{x}}= -\frac{\theta_1 \cdot \sqrt{\pi}}{\sqrt{n}}(\log{2n}+2 \log{2}+C),
\]
while $ 1>\begin{array}{c} \theta \\ ~\theta_1  \end{array}>0$. It indeed becomes clear easily, that it is
\[
\lim{P(n)}=0, ~~~\lim{Q(n)}=0, ~~~[n= \infty],
\]
where it is possible to conclude\\[5mm]
41.
$\begin{cases}
\sin{\frac{1}{3}\pi} \cdot \int_0^{\infty} \frac{\log{x}}{e^x+1+e^{-x}} \cdot \frac{\diff{d}x}{\sqrt{x}} =\sqrt{\pi} \cdot \sum_{i=1}^{\infty} (-1)^i \cdot \sin{\frac{1}{3} i \pi} \cdot \frac{\log{i}+2\log{2}+C}{\sqrt{i}}\\[5mm]
~~~~~~~~~~\int_0^{\infty} \frac{\log{x}}{e^x+e^{-x}} \cdot \frac{\diff{d}x}{\sqrt{x}}=\sqrt{\pi} \cdot \sum_{i=0}^{\infty} (-1)^{i+1}  \cdot \frac{\log{(2i+1)}+2\log{2}+C}{\sqrt{2i+1}}
\end{cases} $\\[2mm]

By the means of the formuilas (39.) we can indeed in an very easy task deduce
\[
\sin{\frac{1}{3}\pi} \cdot \int_0^{\infty} \frac{\diff{d}x}{e^x+1+e^{-x}} \cdot \frac{1}{\sqrt{x}} =\sqrt{\pi} \cdot \sum_{i=1}^{\infty} (-1)^{i-1}  \cdot \frac{ \sin{\frac{1}{3} i \pi}}{\sqrt{i}},
\] 
\[
\int_0^{\infty} \frac{\diff{d}x}{e^x+e^{-x}} \cdot \frac{1}{\sqrt{x}}=\sqrt{\pi} \cdot \sum_{i=0}^{\infty} (-1)^{i}  \cdot \frac{(-1)^{i}}{\sqrt{2i+1}},
\]
which formulas along with (41.), after having taken into account (38.), give after certain very easy reductions\\[5mm]
42.
$\begin{cases}
\frac{\log{1}}{\sqrt{1}}-\frac{\log{2}}{\sqrt{2}}+\frac{\log{4}}{\sqrt{4}}-\frac{\log{5}}{\sqrt{5}}+\frac{\log{7}}{\sqrt{7}}-\text{etc.}\\[5mm]
\frac{1}{2}(\frac{1}{2}\pi-C-\log{\frac{8}{3}\pi})\{\frac{1}{\sqrt{1}}-\frac{1}{\sqrt{2}}+\frac{1}{\sqrt{4}}-\frac{1}{\sqrt{5}}+\frac{1}{\sqrt{7}}-\text{etc.} \}\\[5mm]
\frac{\log{1}}{\sqrt{1}}-\frac{\log{3}}{\sqrt{3}}+\frac{\log{5}}{\sqrt{5}}-\frac{\log{7}}{\sqrt{7}}+\frac{\log{9}}{\sqrt{9}}-\text{etc.}\\[5mm]
\frac{1}{2}(\frac{1}{2}\pi-C-\log{2\pi})\{\frac{1}{\sqrt{1}}-\frac{1}{\sqrt{3}}+\frac{1}{\sqrt{5}}-\frac{1}{\sqrt{7}}+\frac{1}{\sqrt{9}}-\text{etc.} \}
\end{cases}$\\[5mm]
which relations also can be written this way:\\[2mm]
43.
$\begin{cases}
\frac{1^{\frac{1}{\sqrt{1}}}\cdot 4^{\frac{1}{\sqrt{4}}} \cdot 7^{\frac{1}{\sqrt{7}}} \cdot 10^{\frac{1}{\sqrt{10}}} \cdot \cdots}{2^{\frac{1}{\sqrt{2}}}\cdot 5^{\frac{1}{\sqrt{5}}} \cdot 8^{\frac{1}{\sqrt{8}}} \cdot 11^{\frac{1}{\sqrt{11}}} \cdot \cdots}= e^{\frac{1}{2}(\frac{1}{2}\pi-C-\log{\frac{8}{3}\pi})A}\\[5mm]
\frac{1^{\frac{1}{\sqrt{1}}}\cdot 5^{\frac{1}{\sqrt{5}}} \cdot 9^{\frac{1}{\sqrt{9}}} \cdot 13^{\frac{1}{\sqrt{13}}} \cdot \cdots}{3^{\frac{1}{\sqrt{3}}}\cdot 7^{\frac{1}{\sqrt{7}}} \cdot 11^{\frac{1}{\sqrt{11}}} \cdot 15^{\frac{1}{\sqrt{15}}} \cdot \cdots}= e^{\frac{1}{2}(\frac{1}{2}\pi-C-\log{2\pi})B}\\
\end{cases}$\\[5mm]

with 
\[
A=\frac{1}{\sqrt{\pi}} \cdot\int_0^{\infty} \frac{\diff{d}x}{e^x+1+e^{-x}} \cdot \frac{1}{\sqrt{x}}= \frac{1}{\sqrt{1}}-\frac{1}{\sqrt{2}}+\frac{1}{\sqrt{4}}-\frac{1}{\sqrt{5}}+\frac{1}{\sqrt{7}}-\text{etc.}
\]
\[
B=\frac{1}{\sqrt{\pi}} \cdot \int_0^{\infty} \frac{\diff{d}x}{e^x+e^{-x}} \cdot \frac{1}{\sqrt{x}}=\frac{1}{\sqrt{1}}-\frac{1}{\sqrt{3}}+\frac{1}{\sqrt{5}}-\frac{1}{\sqrt{7}}+\frac{1}{\sqrt{9}}-\text{etc.}
\]

\section*{~~~~~~~~~~~~~~~~~~~~~~~~~~~~~~~~~§5}

If we expand $\frac{e^{au}-e^{-au}}{e^{\pi u}-e^{-\pi u}}$ into a series, while it is in identical manner
\[
\frac{1}{e^{\pi u}-e^{-\pi u}}= e^{-\pi u} \cdot \sum_{i=0}^{n-1} e^{-2 i \pi u} +\frac{e^{-2 n \pi u}}{e^{\pi u}-e^{-\pi u}}
\]
it is, of course,
\[
\frac{e^{au}-e^{-au}}{e^{\pi u}-e^{-\pi u}}= \sum_{i=0}^{n-1} [e^{-[(2i+1)\pi -a]u}-e^{-[(2i+1)\pi +a]u}]+\frac{e^{-2n \pi u}(e^{au}-e^{-au})}{e^{\pi u}-e^{-\pi u}},
\]
whence
\[
\int_0^{\infty}\frac{e^{au}-e^{-au}}{e^{\pi u}-e^{-\pi u}} \cdot \frac{\diff{d}u}{u^s}= \Gamma{(1-s)} \cdot \sum_{i=0}^{n-1}[\frac{1}{((2i+1)\pi-a)^{1-s}}+\frac{1}{((2i+1)\pi+a)^{1-s}}]+\varphi{(n)},
\]
where we put for the sake of brevity
\[
\varphi{(n)} = \int_0^{\infty}\frac{e^{au}-e^{-au}}{e^{\pi u}-e^{-\pi u}} \cdot \frac{e^{-2n \pi u} \cdot \diff{d}u}{u^s}= M \cdot \int_0^{\infty}\frac{e^{-2n \pi u} \cdot \diff{d}u}{u^s}= \frac{M\Gamma{(1-s)}}{(2n \pi)^{1-s}},
\]
where $M$ is certain finite quantity. Hence it indeed easily becomes clear
\[
\lim{\varphi{(n)}}=0 ~~~~~~~~ [n = \infty]
\]
and therefore
\[
44. ~~~~ \int_0^{\infty}\frac{e^{au}-e^{-au}}{e^{\pi u}-e^{-\pi u}} \cdot \frac{\diff{d}u}{u^s}= \Gamma{(1-s)} \cdot \sum_{i=0}^{\infty}[\frac{1}{((2i+1)\pi-a)^{1-s}}+\frac{1}{((2i+1)\pi+a)^{1-s}}]
\]
But now, while it is in identical manner
\[
44\frac{1}{2}.~~~~ \frac{\sin{a}}{1+2y\cos{a}+y^2}= \sum_{i=1}^{n} (-1)^{i-1} y^{i-1} \sin{ia} +\frac{(-1)^n \cdot y^n(\sin{(n+1)a}+y\sin{na})}{1+2y\cos{a}+y^2}
\]
we will have in identical manner
\[
45.~~~ \int_0^1 \frac{\sin{a}}{1+2y\cos{a}+y^2} \cdot \frac{\diff{d}y}{(\log{\frac{1}{y}})^{s-1}}
\]
\[
= \Gamma{(s)} \cdot \sum_{i=1}^{n} (-1)^{i-1} \frac{\sin{ia}}{i^s} +(-1)^n(W(n)\sin{(n+1)a}+W(n+1)\sin{na}),
\]
where, for the sake of brevity, we put 
\[
W(n)=\int_0^1 \frac{y^n (\log{\frac{1}{y}})^{s-1}\diff{d}y}{1+2y\cos{a}+y^2}=\theta \cdot \int_0^1 y^n (\log{\frac{1}{y}})^{s-1}\diff{d}y= \frac{\theta \cdot \Gamma{(s)}}{(n+1)^s},
\]
while $1>\theta>0$. It therefore easily becomes clear, that it is
\[
\lim{W(n)} ~~~~~~ [n = \infty];
\] 
whence from (45.) one will obtain
\[
46.~~~~ \int_0^1 \frac{\sin{a}}{1+2y\cos{a}+y^2} \cdot \frac{\diff{d}y}{(\log{\frac{1}{y}})^{s-1}} = \Gamma{(s)}\cdot \sum_{i=1}^{\infty} (-1)^{i-1} \cdot \frac{\sin{ia}}{i^s}
\]

But after having substituted the value in (29.), which the formulas (44. and 46.) give, we have this remarkable relation between two infinte series, if $s$ is changed into $1-s$:\\

47.
$\begin{cases}
\frac{1}{(\pi-a)^s}-\frac{1}{(\pi+a)^s}+\frac{1}{(3\pi-a)^s}-\frac{1}{(3\pi+a)^s}+\frac{1}{(5\pi-a)^s}-\frac{1}{(5\pi+a)^s}+\text{etc.}\\[3mm]
\frac{1}{\sin{\frac{1}{2} s \pi} \cdot \Gamma{(s)}}\{ \frac{\sin{a}}{1^{1-s}}-\frac{\sin{2a}}{2^{1-s}}+\frac{\sin{3a}}{3^{1-s}}-\frac{\sin{4a}}{4^{1-s}}+\frac{\sin{5a}}{5^{1-s}}-\text{etc.}\}
\end{cases}$\\[5mm]

and if we put $\pi-a$ in the place of $a$:\\

48.
$\begin{cases}
\frac{1}{a^s}-\frac{1}{(2\pi-a)^s}+\frac{1}{(2\pi+a)^s}-\frac{1}{(4\pi-a)^s}+\frac{1}{(4\pi+a)^s}-\frac{1}{(6\pi-a)^s}+\text{etc.}\\[3mm]
\frac{1}{\sin{\frac{1}{2} s \pi} \cdot \Gamma{(s)}}\{ \frac{\sin{a}}{1^{1-s}}+\frac{\sin{2a}}{2^{1-s}}+\frac{\sin{3a}}{3^{1-s}}+\frac{\sin{4a}}{4^{1-s}}+\frac{\sin{5a}}{5^{1-s}}+\text{etc.}\}
\end{cases}$\\[5mm]

Hence, if we set $s=\frac{1}{2}$, it is of course\\

49.
$\begin{cases}
\frac{1}{\sqrt{\pi-a}}-\frac{1}{\sqrt{\pi+a}}+\frac{1}{\sqrt{3\pi-a}}-\frac{1}{\sqrt{3\pi+a}}+\frac{1}{\sqrt{5\pi-a}}-\frac{1}{\sqrt{5\pi+a}}+\text{etc.}\\[3mm]
=\sqrt{\frac{2}{\pi}}\{ \frac{\sin{a}}{\sqrt{1}}-\frac{\sin{2a}}{\sqrt{2}}+\frac{\sin{3a}}{\sqrt{3}}-\frac{\sin{4a}}{\sqrt{4}}+\frac{\sin{5a}}{\sqrt{5}}-\text{etc.}\}\\[3mm]
\frac{1}{\sqrt{a}}-\frac{1}{\sqrt{2\pi-a}}+\frac{1}{\sqrt{2\pi+a}}-\frac{1}{\sqrt{4\pi-a}}+\frac{1}{\sqrt{4\pi+a}}-\frac{1}{\sqrt{6\pi-a}}+\text{etc.}\\[3mm]
=\sqrt{\frac{2}{\pi}}\{ \frac{\sin{a}}{\sqrt{1}}+\frac{\sin{2a}}{\sqrt{2}}+\frac{\sin{3a}}{\sqrt{3}}+\frac{\sin{4a}}{\sqrt{4}}+\frac{\sin{5a}}{\sqrt{5}}+\text{etc.}\}
\end{cases}$\\[2mm]

Let us suppose in (47. and 48.), that $a$ is in a commensurable ratio to $\pi$, i.e. $a=\frac{m \pi}{n}$ $(m<n ~~\text{integer number})$; then one will easily obtain\\

50.
$\begin{cases}
\frac{1}{(n-m)^s}-\frac{1}{(n+m)^s}+\frac{1}{(3n+m)^s}-\frac{1}{(3n+m)^s}+\frac{1}{(5n-m)^s}-\frac{1}{(5n+m)^s}+\text{etc.}\\[3mm]
=\frac{(\frac{\pi}{n})^s}{\sin{\frac{1}{2} s \pi} \cdot \Gamma{(s)}} \{ \frac{\sin{\frac{m \pi}{n}}}{1^{1-s}}-\frac{\sin{\frac{2m \pi}{n}}}{2^{1-s}}+\frac{\sin{\frac{3m \pi}{n}}}{3^{1-s}}-\frac{\sin{\frac{4m \pi}{n}}}{4^{1-s}}+\text{etc.}\}\\[3mm]
\frac{1}{m^s}-\frac{1}{(2n-m)^s}+\frac{1}{(2n+m)^s}-\frac{1}{(4n-m)^s}+\frac{1}{(4n+m)^s}-\frac{1}{(6n-m)^s}+\text{etc.}\\[3mm]
=\frac{(\frac{\pi}{n})^s}{\sin{\frac{1}{2} s \pi} \cdot \Gamma{(s)}} \{ \frac{\sin{\frac{m \pi}{n}}}{1^{1-s}}+\frac{\sin{\frac{2m \pi}{n}}}{2^{1-s}}+\frac{\sin{\frac{3m \pi}{n}}}{3^{1-s}}+\frac{\sin{\frac{4m \pi}{n}}}{4^{1-s}}+\text{etc.}\}
\end{cases} $\\[5mm]

Ex. 1. While $m=1$, $n=2$, it is
\[
51. ~~~ \frac{1}{1^s}-\frac{1}{3^s}+\frac{1}{5^s}-\frac{1}{7^s}+\frac{1}{9^s}-\text{etc.} \cdots
\]
\[
=\frac{(\frac{1}{2}\pi)^s}{\sin{\frac{1}{2} s \pi} \cdot \Gamma{(s)}}\{ \frac{1}{1^{1-s}}-\frac{1}{3^{1-s}}+\frac{1}{5^{1-s}}-\frac{1}{7^{1-s}}+\text{etc.}\}
\]
Ex. 2. While in the first of the formulas it is $m=1$, $n=3$, it is
\[
52. ~~~ \frac{1}{1^s}-\frac{1}{2^s}+\frac{1}{4^s}-\frac{1}{5^s}+\frac{1}{7^s}-\text{etc.} \cdots
\]
\[
=\frac{(\frac{2}{3}\pi)^s \cdot \sin{\frac{1}{3}\pi}}{\sin{\frac{1}{2} s \pi} \cdot \Gamma{(s)}}\{ \frac{1}{1^{1-s}}-\frac{1}{2^{1-s}}+\frac{1}{4^{1-s}}-\frac{1}{5^{1-s}}+\frac{1}{7^{1-s}}-\text{etc.}\}
\]\\

I remember, that I saw the formulas (51. and 52.) (if I recall correctly\footnote{Malmstèn is referring to Euler's paper \cite{1}, but Euler found the functional equation for the function $\eta{(s)}= \sum_{i=1}^{\infty} \frac{(-1)^{i-1}}{i^s} $ and the one for Malmstén's formula (51.)). Euler gave the functional equation for $\eta{(s)}$ in \cite{2} again. }), found elsewhere by Euler by induction, but free from any proof; and we do not find another proof of them at anyone, although in their form the seem to be worth of the Geometers' attention.\\

Ex. 3. While in the second of the formulas (50.) it is $m=1$, $n=3$, after having put for the sake of brevity
\begin{align*}
f(s)&= \frac{1}{1^s}-\frac{1}{5^s}+\frac{1}{7^s}-\frac{1}{11^s}+\frac{1}{13^s}-\frac{1}{17^s}-\frac{1}{19^s}-\text{etc.}\\
\varphi{(s)}&= \frac{1}{1^s}+\frac{1}{2^s}-\frac{1}{4^s}-\frac{1}{5^s}+\frac{1}{7^s}+\frac{1}{8^s}-\frac{1}{10^s}-\frac{1}{11^s}+\text{etc.}
\end{align*}
it is, of course,
\[
f(s)= \frac{(\frac{2}{3}\pi)^s \cdot \sin{\frac{1}{3}\pi}}{\sin{\frac{1}{2} s \pi} \cdot \Gamma{(s)}} \cdot \varphi{(1-s)}.
\]
But because it is
\[
\varphi{(s)}= f(s)+2^{-s}\{ \frac{1}{1^s}- \frac{1}{2^s}+ \frac{1}{4^s}- \frac{1}{5^s}+\text{etc.}\}
\]
and also
\[
\varphi{(s)}= f(s)+2^{-s}f(s)-2^{-2s}\{ \frac{1}{1^s}- \frac{1}{2^s}+ \frac{1}{4^s}- \frac{1}{5^s}+\text{etc.}\},
\]
whence it easily emerges
\[
(1+2^{-s})\varphi{(s)}= (1+2^{1-s})f(s),
\]
we will have these two relations\\

53.
$\begin{cases}
f(s)=\frac{(\frac{1}{3}\pi)^s \cdot \sin{\frac{1}{3}\pi}}{\sin{\frac{1}{2} s \pi} \cdot \Gamma{(s)}} \cdot \frac{1+2^s}{1+2^{s-1}} \cdot f(1-s)\\[3mm]
\varphi{(s)}=\frac{(\frac{1}{3}\pi)^s \cdot \sin{\frac{1}{3}\pi}}{\sin{\frac{1}{2} s \pi} \cdot \Gamma{(s)}} \cdot \frac{2(1+2^{s-1})}{1+2^{s}} \cdot \varphi{(1-s)}
\end{cases}$\\[5mm]

Ex. 4. While in (50.) it is $m=1$, $m=4$, after having put for the sake of brevity:
\begin{align*}
F(s)&= \frac{1}{1^s}+\frac{1}{3^s}-\frac{1}{5^s}-\frac{1}{7^s}+\frac{1}{9^s}+\frac{1}{11^s}-\frac{1}{13^s}-\frac{1}{15^s}+\text{etc.}\\
W(s)&= \frac{1}{1^s}-\frac{1}{3^s}+\frac{1}{5^s}-\frac{1}{7^s}+\frac{1}{9^s}-\frac{1}{11^s}+\text{etc.}\\
\psi{(s)}&= \frac{1}{3^s}-\frac{1}{5^s}+\frac{1}{11^s}-\frac{1}{13^s}+\frac{1}{19^s}-\frac{1}{21^s}+\text{etc.}\\
P(s)&= \frac{1}{1^s}-\frac{1}{7^s}+\frac{1}{9^s}-\frac{1}{15^s}+\frac{1}{17^s}-\text{etc.}
\end{align*}\\

it is , of course
\begin{align*}
\psi{(s)}= \frac{(\frac{\pi}{4})^s \cdot \sin{\frac{\pi}{4}}}{\sin{\frac{1}{2} s \pi} \cdot \Gamma{(s)}}\cdot F(1-s)-\frac{\frac{1}{2} \cdot (\frac{1}{2}\pi)^s}{\sin{\frac{1}{2} s \pi} \cdot \Gamma{(s)}}\cdot W(1-s)\\
P(s)=\frac{(\frac{\pi}{4})^s \cdot \sin{\frac{\pi}{4}}}{\sin{\frac{1}{2} s \pi} \cdot \Gamma{(s)}}\cdot W(1-s)+\frac{\frac{1}{2} \cdot (\frac{1}{2}\pi)^s}{\sin{\frac{1}{2} s \pi} \cdot \Gamma{(s)}}\cdot F(1-s)
\end{align*}\\

whence by adding, because it is 
\[
\psi{(s)}+P(s)= F(s),
\]
it will be 
\[
54. ~~~ F(s)= \frac{2 \cdot (\frac{\pi}{4})^s \cdot \sin{\frac{\pi}{4}}}{\sin{\frac{1}{2} s \pi} \cdot \Gamma{(s)}}\cdot F(1-s)
\]\\

\section*{~~~~~~~~~~~~~~~~~~~~~~~~~~~~~~~~~§6}

Let us differentiate formula (44.) with respect to $s$ as variable; then it will be
\[
\int_0^{\infty}\frac{e^{au}-e^{-au}}{e^{\pi u}-e^{-\pi u}} \cdot u^{-s} \log{u} \cdot \diff{d}u
\]
\[
=Z^{\prime}(1-s) \cdot \int_0^{\infty}\frac{e^{au}-e^{-au}}{e^{\pi u}-e^{-\pi u}} \cdot \frac{\diff{d}u}{u^s}-\Gamma{(1-s)}\cdot \sum_{i=0}^{\infty} \bigg[\frac{\log{((2i+1)\pi-a)}}{((2i+1)\pi-a)^{1-s}}-\frac{\log{((2i+1)\pi+a)}}{((2i+1)\pi+a)^{1-s}}\bigg],
\]\\
whence for $s=0$, while
\[
Z^{\prime}(1)= -C~~~ \text{and} ~~~ \int_0^{\infty}\frac{e^{au}-e^{-au}}{e^{\pi u}-e^{-\pi u}} \cdot \diff{d}u= \tan{\frac{1}{2}a},
\]\\
we will have
\[
\sum_{i=0}^{\infty} \bigg[\frac{\log{((2i+1)\pi-a})}{((2i+1)\pi-a)^{1-s}}-\frac{\log{((2i+1)\pi+a)}}{((2i+1)\pi+a)^{1-s}}\bigg]
\]
\[
=-\frac{1}{2}C \cdot \tan{\frac{1}{2}a} - \int_0^{\infty}\frac{e^{au}-e^{-au}}{e^{\pi u}-e^{-\pi u}}\log{u} \cdot \diff{d}u.
\]\\

But in whatever commensurable ratio $a$ is to $\pi$, the formulas (7.) (for $x=0$) give the value for the integral, that occurs on the right hand side. Therefore let be $a=\frac{m \pi}{n}$ $m$ and $n$ integer numbers $n>n$);\\
\[
\frac{n}{\pi} \cdot \sum_{i=0}^{\infty} \bigg[\frac{\log{((2i+1)n-m)}}{((2i+1)n-m)^{1-s}}-\frac{\log{((2i+1)n+m)}}{((2i+1)n+m)^{1-s}}\bigg]
\]
\[
\log{\frac{n}{\pi}} \cdot \sum_{i=0}^{\infty} \bigg[\frac{1}{((2i+1)\pi-\frac{m \pi}{n})^{1-s}}-\frac{1}{((2i+1)\pi+\frac{m \pi}{n})^{1-s}}\bigg]
\]
\[
=-\frac{1}{2}C \cdot \tan{\frac{m \pi}{2n}} - \int_0^{\infty}\frac{e^{\frac{m \pi}{n} \cdot u}-e^{-\frac{m \pi}{n} \cdot u}}{e^{\pi u}-e^{-\pi u}}\cdot \log{u} \cdot \diff{d}u.
\]\\
whence, because it is from (44.) (for $s=0$)
\[
\sum_{i=0}^{\infty} \bigg[\frac{1}{((2i+1)\pi-\frac{m \pi}{n})^{1-s}}-\frac{1}{((2i+1)\pi+\frac{m \pi}{n})^{1-s}}\bigg]= \int_0^{\infty}\frac{e^{\frac{m \pi}{n} \cdot u}-e^{-\frac{m \pi}{n} \cdot u}}{e^{\pi u}-e^{-\pi u}} \cdot \diff{d}u= \frac{1}{2}\tan{\frac{1}{2}a},
\]\\

it finally is from the cited formulas (7.)\\

55.
$\begin{cases}
\frac{\log{(n-m)}}{n-m}-\frac{\log{(n+m)}}{n+m}+\frac{\log{(3n-m)}}{3n-m}-\frac{\log{(3n+m)}}{3n+m}+\frac{\log{(5n-m)}}{5n-m}-\text{etc.}\\[3mm]
= -\frac{\pi}{2n} \cdot \tan{\frac{m \pi}{2n}}(C+\log{2 \pi})- \frac{\pi}{n} \cdot \sum_{i=1}^{n-1} (-1)^{i-1} \sin{\frac{im \pi}{n}}\log\bigg\{{\frac{\Gamma{(\frac{n+i}{2n}})}{\Gamma{(\frac{i}{2n}})}}\bigg\}\\[5mm]
~~~~~~~~~~~~~~~~~~~~~~~~~(m+n= \text{odd number})\\[5mm]
\frac{\log{(n-m)}}{n-m}-\frac{\log{(n+m)}}{n+m}+\frac{\log{(3n-m)}}{3n-m}-\frac{\log{(3n+m)}}{3n+m}+\frac{\log{(5n-m)}}{5n-m}-\text{etc.}\\[3mm]
= -\frac{\pi}{2n} \cdot \tan{\frac{m \pi}{2n}}(C+\log{ \pi})- \frac{\pi}{n} \cdot \sum_{i=1}^{\frac{1}{2}(n-1)} (-1)^{i-1} \sin{\frac{im \pi}{n}}\log\bigg\{{\frac{\Gamma{(\frac{n-i}{n}})}{\Gamma{(\frac{i}{n}})}}\bigg\}\\[5mm]
~~~~~~~~~~~~~~~~~~~~~~~~~(m+n= \text{even number}),\\
\end{cases}$\\[5mm]
and if we put $n-m$ in the place of $m$,\\[5mm]

56.
$\begin{cases}
\frac{\log{m}}{m}-\frac{\log{(2n-m)}}{2n-m}+\frac{\log{(2n+m)}}{2n+m}-\frac{\log{(4n-m)}}{4n-m}+\frac{\log{(4n+m)}}{4n+m}-\text{etc.}\\[3mm]
= -\frac{\pi}{2n} \cdot \cot{\frac{m \pi}{2n}}(C+\log{2 \pi})- \frac{\pi}{n} \cdot \sum_{i=1}^{n-1}  \sin{\frac{im \pi}{n}}\log\bigg\{{\frac{\Gamma{(\frac{n+i}{2n}})}{\Gamma{(\frac{i}{2n}})}}\bigg\}\\[5mm]
~~~~~~~~~~~~~~~~~~~~~~~~~(m= \text{even number})\\[5mm]
\frac{\log{m}}{m}-\frac{\log{(2n-m)}}{2n-m}+\frac{\log{(2n+m)}}{2n+m}-\frac{\log{(4n-m)}}{4n-m}+\frac{\log{(4n+m)}}{4n+m}-\text{etc.}\\[3mm]
= -\frac{\pi}{2n} \cdot \cot{\frac{m \pi}{2n}}(C+\log{ \pi})- \frac{\pi}{n} \cdot \sum_{i=1}^{\frac{1}{2}(n-1)}  \sin{\frac{im \pi}{n}}\log\bigg\{{\frac{\Gamma{(\frac{n-i}{n}})}{\Gamma{(\frac{i}{n}})}}\bigg\}\\[5mm]
~~~~~~~~~~~~~~~~~~~~~~~~~(m= \text{odd number}).\\
\end{cases}$\\[5mm]

Ex. 1. After having put $n=2$, $m=1$, it is
\[
\frac{\log{1}}{{1}}-\frac{\log{3}}{{3}}+\frac{\log{5}}{{5}}-\frac{\log{7}}{{7}}+\text{etc.}= \frac{\pi}{4}(\log{\pi}-C)-\pi \log{\Gamma{(\frac{3}{4})}}
\]
and therefore
\[
\frac{{1}\cdot 5^{\frac{1}{{5}}} \cdot 9^{\frac{1}{9}} \cdot 13^{\frac{1}{{13}}} \cdot \cdots}{3^{\frac{1}{{3}}}\cdot 7^{\frac{1}{{7}}} \cdot 11^{\frac{1}{{11}}} \cdot 15^{\frac{1}{{15}}} \cdot \cdots}=\bigg\{\frac{\pi e^{-C}}{(\Gamma{(\frac{3}{4})})^4}\bigg\}^{\frac{\pi}{4}}
\]\\
Ex. 2. After having put $n=3$, $m=1$, while
\[
1-\frac{1}{2}+\frac{1}{4}-\frac{1}{5}+\frac{1}{7}-\frac{1}{8}+\text{etc.}= \frac{\pi}{3\sqrt{3}}
\]
it is of course
\[
\frac{\log{1}}{{1}}-\frac{\log{2}}{{2}}+\frac{\log{4}}{{4}}-\frac{\log{5}}{{5}}+\frac{\log{7}}{{7}}-\text{etc.}=\frac{\pi}{\sqrt{3}}\log\bigg\{{\frac{\Gamma{(\frac{1}{3})}}{\Gamma{(\frac{2}{3})}}-\frac{\pi}{3\sqrt{3}}(C+\log{2 \pi}})\bigg\},
\]
\[
\frac{\log{1}}{{1}}-\frac{\log{5}}{{5}}+\frac{\log{7}}{{7}}-\frac{\log{11}}{{11}}+\frac{\log{13}}{{13}}-\text{etc.}
\]
\[
=\frac{\pi}{2\sqrt{3}}\bigg\{ \log{\frac{2 \pi}{\sqrt{3}}}-C-2\log[\Gamma{(\frac{5}{6})}\cdot \Gamma{(\frac{2}{3})}]\bigg\}
\]\\

and from this 
\begin{align*}
\frac{{1}\cdot 4^{\frac{1}{{4}}} \cdot 7^{\frac{1}{7}} \cdot 10^{\frac{1}{{10}}} \cdot \cdots}{2^{\frac{1}{{2}}}\cdot 5^{\frac{1}{{5}}} \cdot 8^{\frac{1}{{8}}} \cdot 11^{\frac{1}{{11}}} \cdot \cdots}&=\bigg\{\frac{\Gamma{(\frac{1}{3})} e^{-\frac{1}{3}C}}{\Gamma{(\frac{2}{3})}\cdot (2 \pi)^{\frac{1}{3}}}\bigg\}^{\frac{\pi}{\sqrt{3}}},\\
\frac{{1}\cdot 7^{\frac{1}{{7}}} \cdot 13^{\frac{1}{13}} \cdot 19^{\frac{1}{{19}}} \cdot \cdots}{5^{\frac{1}{{5}}}\cdot 11^{\frac{1}{{11}}} \cdot 17^{\frac{1}{{17}}} \cdot 23^{\frac{1}{{23}}} \cdot \cdots}&=\bigg\{\frac{ (2 \pi)^{\frac{1}{2}} \cdot e^{-\frac{1}{2}C}}{3^{\frac{1}{4}} \cdot \Gamma{(\frac{2}{3})}\cdot \Gamma{(\frac{5}{6})}} \bigg\}^{\frac{\pi}{\sqrt{3}}}.
\end{align*}\\

Ex. 3. After having put $n=4$, $m=1$, it is
\begin{align*}
\frac{\log{1}}{{1}}-\frac{\log{3}}{{3}}+\frac{\log{5}}{{5}}-\frac{\log{11}}{{11}}+\frac{\log{13}}{{13}}-\text{etc.}&= \frac{\pi}{2\sqrt{2}}\log\bigg\{{\frac{\Gamma{(\frac{5}{8})}\cdot \Gamma{(\frac{3}{8})}}{2^{\frac{1}{4}}\cdot \pi^{\frac{1}{2}} \cdot e^{\frac{1}{3}C}}}\bigg\}+\frac{1}{2}\pi\log\bigg\{{\frac{e^{-\frac{1}{4}C}\cdot \pi^{\frac{1}{4}}}{\Gamma{(\frac{3}{4})}}}\bigg\},\\
\frac{\log{1}}{{1}}-\frac{\log{7}}{{7}}+\frac{\log{9}}{{9}}-\frac{\log{15}}{{15}}+\frac{\log{17}}{{17}}-\text{etc.}&= \frac{1}{2}\pi\log\bigg\{{\frac{\pi^{\frac{1}{4}}e^{-\frac{1}{4}}C}{\Gamma{(\frac{3}{4})}}}\bigg\}+\frac{\pi}{2\sqrt{2}}\log\bigg\{{\frac{2^{\frac{1}{4}}\cdot \pi^{\frac{1}{2}} \cdot e^{\frac{1}{3}C}}{\Gamma{(\frac{5}{8})}\cdot \Gamma{(\frac{3}{8})}}}\bigg\}
\end{align*}\\

and from this
\begin{align*}
\frac{{1}\cdot 5^{\frac{1}{{5}}} \cdot 13^{\frac{1}{13}} \cdot 21^{\frac{1}{{21}}} \cdot \cdots}{3^{\frac{1}{{3}}}\cdot 11^{\frac{1}{{11}}} \cdot 19^{\frac{1}{{19}}} \cdot 27^{\frac{1}{{27}}} \cdot \cdots}=\bigg(\frac{e^{-\frac{1}{4}C}\cdot \pi^{\frac{1}{4}}}{\Gamma{(\frac{3}{4})}}\bigg)^{\frac{1}{2}\pi} \cdot \bigg\{ \frac{\Gamma{(\frac{5}{8})}\cdot \Gamma{(\frac{3}{8})}}{2^{\frac{1}{4}}\cdot \pi^{\frac{1}{2}} \cdot e^{\frac{1}{3}C}}\bigg\}^{\frac{\pi}{2\sqrt{2}}},\\
\frac{{1}\cdot 9^{\frac{1}{{9}}} \cdot 17^{\frac{1}{17}} \cdot 25^{\frac{1}{{25}}} \cdot \cdots}{7^{\frac{1}{{7}}}\cdot 15^{\frac{1}{{15}}} \cdot 23^{\frac{1}{{23}}} \cdot 31^{\frac{1}{{31}}} \cdot \cdots}=\bigg(\frac{e^{-\frac{1}{4}C}\cdot \pi^{\frac{1}{4}}}{\Gamma{(\frac{3}{4})}}\bigg)^{\frac{1}{2}\pi} \cdot \bigg\{ \frac{2^{\frac{1}{4}}\cdot \pi^{\frac{1}{2}} \cdot e^{\frac{1}{3}C}}{\Gamma{(\frac{5}{8})}\cdot \Gamma{(\frac{3}{8})}}\bigg\}^{\frac{\pi}{2\sqrt{2}}}.
\end{align*}\\

\section*{~~~~~~~~~~~~~~~~~~~~~~~~~~~~~~~~~§7}

If one puts $a=0$ in formula (29.), it is
\[
57.~~~ \int_0^{\infty} \frac{2u \cdot \diff{d}u}{e^{\pi u}-e^{-\pi u}}\cdot \frac{\cos{(s\arctan{\frac{u}{x}})}}{(x^2+u^2)^{\frac{1}{2}s}} = \frac{1}{\Gamma{(s)}}\cdot \int_0^{1} \frac{y^x(\log\frac{1}{y})^{s-1}\diff{d}y}{(1+y)^2},
\]
whence we will easily have by subtraction
\[
58. ~~~ \int_0^{\infty} \frac{e^{au}-e^{-au}-2au} {e^{\pi u}-e^{-\pi u}}\cdot  \frac{\cos{(s\arctan{\frac{u}{x}})}}{(x^2+u^2)^{\frac{1}{2}s}}\cdot \diff{d}u
\]
\[
= \frac{1}{\Gamma{(s)}}\cdot \int_0^{1} \bigg(\frac{\sin{a}}{1+2y\cos{a}+y^2}-\frac{a}{(1+y)^2}\bigg) \cdot y^x(\log\frac{1}{y})^{s-1}\diff{d}y.
\]
Let us differentiate this formula with respect to $s$ as a variable, and let us put $s=1$ afterwards; then, because 
\[
\arctan{\frac{u}{x}}= \frac{1}{2}\pi-\arctan{\frac{x}{u}}, ~~~ Z^{\prime}(1)=-C,
\]
we will obtain\\

59.
$\begin{cases}
\frac{1}{2}\pi \int_0^{\infty} \frac{2u \cdot \diff{d}u}{e^{\pi u}-e^{-\pi u}}\cdot \frac{u \diff{d}u}{x^2+u^2}-xF(x)+\frac{1}{2}W(x)\\[3mm]
= -\int_0^{1} \bigg(\frac{\sin{a}}{1+2y\cos{a}+y^2}-\frac{a}{(1+y)^2}\bigg) \cdot y^x\log(\log\frac{1}{y})\diff{d}y\\[3mm]
~~~~~+C\int_0^{1} \bigg(\frac{\sin{a}}{1+2y\cos{a}+y^2}-\frac{a}{(1+y)^2}\bigg) \cdot y^x\cdot\diff{d}y
\end{cases}$\\[3mm]

after having put for the sake of brevity
\[
F(x)=  \int_0^{\infty} \frac{e^{au}-e^{-au}-2au}{e^{\pi u}-e^{-\pi u}} \cdot \frac{\arctan{\frac{u}{x}}}{\frac{u}{x}}\cdot \frac{ \diff{d}u}{x^2+u^2}
\]
\[
\theta  \int_0^{\infty} \frac{e^{au}-e^{-au}-2au}{e^{\pi u}-e^{-\pi u}} \cdot \frac{ \diff{d}u}{x^2+u^2},~~~ (1>\theta>0)
\]
\[
W(x)= \int_0^{\infty} \frac{e^{au}-e^{-au}-2au}{e^{\pi u}-e^{-\pi u}}\cdot \frac{\log{(x^2+u^2} \diff{d}u}{x^2+u^2}.
\]
But because it is
\[
F(x)= \frac{1}{2}\theta \int_0^{\infty} \frac{e^{au}-e^{-au}-2au}{u(e^{\pi u}-e^{-\pi u})} \cdot \frac{2u \diff{d}u}{x^2+u^2}+\theta \int_0^{\infty} \frac{e^{au}-e^{-au}-2au}{e^{\pi u}-e^{-\pi u}} \cdot \frac{ \diff{d}u}{x^2+u^2}
\]
\[
=M\int_0^{1}\frac{2u \diff{d}u}{x^2+u^2}+ \frac{\theta}{x^2+\xi^2} \cdot \int_0^{1} \frac{e^{au}-e^{-au}-2au}{u(e^{\pi u}-e^{-\pi u})} \cdot \diff{d}u
\]
\[
=M(\log{(1+x^2)}-\log{x^2})+\frac{\theta N}{x^2+\xi^2}
\]
(where, whatever $x$ is, $M$ and $N$ always retain a finite value and $\xi$ is a quantity $<1$), and
\[
W(x)= \frac{1}{2} \int_0^{\infty} \frac{e^{au}-e^{-au}-2au}{u(e^{\pi u}-e^{-\pi u})} \cdot \frac{\log{(x^2+u^2)} \cdot 2u \diff{d}u}{x^2+u^2}+ \int_0^{\infty} \frac{e^{au}-e^{-au}-2au}{e^{\pi u}-e^{-\pi u}} \cdot \frac{\log{(x^2+u^2)} \diff{d}u}{x^2+u^2}
\]
\[
=M_1\int_0^{1}\frac{\log{(x^2+u^2)} \cdot 2u \diff{d}u}{x^2+u^2}+ \frac{\log{(x^2+\xi_1^2)}}{x^2+\xi_1^2} \cdot \int_0^{1} \frac{e^{au}-e^{-au}-2au}{u(e^{\pi u}-e^{-\pi u})} \cdot \diff{d}u
\]
\[
=M_1((\log{(1+x^2)})^2-(\log{x^2})^2)+N \cdot\frac{\log{(x^2+\xi_1^2)}}{x^2+\xi_1^2}
\]\\

(where, whatever $x$ is, completely the same holds about $M_1$ and $\xi_1$, what was said about $M$ and $\xi$ above); it is easy to conclude
\[
\lim{xF(x)}=0,~~~\lim{xW(x)}=0, ~~~~[x=0],
\]
whence we will have from (59.), while $x$ converges to zero,
\[
\frac{1}{2}\pi\int_0^{\infty} \frac{e^{au}-e^{-au}-2au}{e^{\pi u}-e^{-\pi u}} \cdot \frac{\diff{d}u}{u}
\]
\[
= -\int_0^{1} \bigg(\frac{\sin{a}}{1+2y\cos{a}+y^2}-\frac{a}{(1+y)^2}\bigg) \log(\log\frac{1}{y})\diff{d}y~+C \cdot\int_0^{1} \bigg(\frac{\sin{a}}{1+2y\cos{a}+y^2}-\frac{a}{(1+y)^2}\bigg)  \diff{d}y
\]
or
\[
60.~~~\int_0^{\infty} \frac{e^{au}-e^{-au}-2au}{e^{\pi u}-e^{-\pi u}} \cdot \frac{\diff{d}u}{u}=-\frac{2}{\pi}\int_0^{1} \bigg(\frac{\sin{a}}{1+2y\cos{a}+y^2}-\frac{a}{(1+y)^2}\bigg) \log(\log\frac{1}{y})\diff{d}y,
\]
because
\[
\int_0^{1} \bigg(\frac{\sin{a}}{1+2y\cos{a}+y^2}-\frac{a}{(1+y)^2}\bigg)  \diff{d}y=0
\]
Let us set
\[
K=\int_0^{\infty} \frac{e^{au}-e^{-au}-2au}{e^{\pi u}-e^{-\pi u}} \cdot \frac{\diff{d}u}{u};
\]
after having differentiated with respect to $a$, we will have 
\[
\frac{\diff{d}K}{\diff{d}a}= \int_0^{\infty} \frac{e^{au}+e^{-au}-2}{e^{\pi u}-e^{-\pi u}} \cdot \diff{d}u=\frac{1}{\pi} \int_0^1 \frac{y^{\frac{a}{\pi}}+y^{-\frac{a}{\pi}}-2}{1-y^2}\diff{d}y,
\]
whence
\[
61. ~~~ \frac{\diff{d}K}{\diff{d}a}= \frac{1}{2\pi}\bigg\{ 2Z^{\prime}(\frac{1}{2})-Z^{\prime}(\frac{1}{2}+\frac{a}{2 \pi})-Z^{\prime}(\frac{1}{2}-\frac{a}{2 \pi})\bigg\},
\]
because it is in general \footnote{See Legendre Exerc. du Calc. Int, book 2, page 156, C.J. Malmstén} 
\[
\int_0^1 \frac{y^b-y^b}{1-y^2}= \frac{1}{2}Z^{\prime}(\frac{1}{2}(b+1))-\frac{1}{2}Z^{\prime}(\frac{1}{2}(a+1))
\] 
From formula (61.), after having integrated from $a=0$, it indeed flows
\[
\int_0^{\infty} \frac{e^{au}-e^{-au}-2au}{e^{\pi u}-e^{-\pi u}} \cdot \frac{\diff{d}u}{u}=\frac{a}{\pi} \cdot Z^{\prime}(\frac{1}{2})-\log\bigg\{{\frac{\Gamma({\frac{1}{2}+\frac{a}{2 \pi}})}{\Gamma({\frac{1}{2}-\frac{a}{2 \pi}})}}\bigg\},
\]
which, compared to (60.), gives 
\[
62.~~~\int_0^{1} \frac{\sin{a} \cdot \log(\log\frac{1}{y})\diff{d}y}{1+2y\cos{a}+y^2}=a\int_0^{1}  \frac{  \log(\log\frac{1}{y})\diff{d}y}{(1+y)^2}-\frac{1}{2}aZ^{\prime}(\frac{1}{2})+\frac{1}{2}\pi\log\bigg\{{\frac{\Gamma({\frac{1}{2}+\frac{a}{2 \pi}})}{\Gamma({\frac{1}{2}-\frac{a}{2 \pi}})}}\bigg\}
\]
For $a=\frac{1}{2}\pi$ by means of formula (12.) (the first) we will obtain
\[
\int_0^{1}  \frac{  \log(\log\frac{1}{y})\diff{d}y}{(1+y)^2} =\frac{1}{2}Z^{\prime}(\frac{1}{2})+\frac{1}{2}\log{2\pi},
\]
which, substituted in (62.), finally gives
\[
63.~~~ \int_0^{1} \frac{  \log(\log\frac{1}{y})\diff{d}y}{1+2y\cos{a}+y^2}= \frac{\pi}{\sin{a}} \cdot \log\bigg\{{\frac{(2\pi)^{\frac{a}{\pi}}\cdot \Gamma({\frac{1}{2}+\frac{a}{2 \pi}})}{\Gamma({\frac{1}{2}-\frac{a}{2 \pi}})}}\bigg\};
\]
this new formula seems quite remarkable.

\section*{~~~~~~~~~~~~~~~~~~~~~~~~~~~~~~~~~§8}

If formula (44$\frac{1}{2}$.) is multiplied by $\log(\log\frac{1}{y})$, after having integrated from $y=0$ to $y=1$, we will have from formula (63.)
\[
\sum_{i=1}^n (-1)^{i-1} \sin{ia} \int_0^1 y^{i-1} \log(\log\frac{1}{y}) \diff{d}y
\]
\[
+(-1)^n\{\sin{(n+1)a} \cdot G(n)+\sin{na} \cdot G(n+1)\}
\]
\[
=\frac{1}{2}\pi\log\bigg\{{\frac{(2\pi)^{\frac{a}{\pi}}\cdot \Gamma({\frac{1}{2}+\frac{a}{2 \pi}})}{\Gamma({\frac{1}{2}-\frac{a}{2 \pi}})}}\bigg\},
\]
after having put for the sake of brevity
\[
G(n)=\int_0^{1} \frac{ y^n \log(\log\frac{1}{y})\diff{d}y}{1+2y\cos{a}+y^2}=\theta\cdot \int_0^{1}y^n \log(\log\frac{1}{y})\diff{d}y,~~~(1>\theta>0),
\]
whence, because
\[
\int_0^{1}y^{r-1} \log(\log\frac{1}{y})\diff{d}y=-\frac{\log{r}+C}{r}
\]
(where $C$ is the Eulerian constant $0,577216\cdots$), we will obtain 
\[
\sum_{i=1}^n (-1)^{i-1} \cdot \frac{\sin{ia}\log{i}}{i}+C \cdot \sum_{i=1}^n (-1)^{i-1}  \frac{\sin{ia}}{i}
\]
\[
+(-1)^{n+1}\{\sin{(n+1)a} \cdot G(n)+\sin{na} \cdot G(n+1)\}
\]
\[
=\frac{1}{2}\pi\log\bigg\{{\frac{ \Gamma({\frac{1}{2}+\frac{a}{2 \pi}})}{\Gamma({\frac{1}{2}-\frac{a}{2 \pi}})}}\bigg\}+\frac{1}{2}a\log{2\pi}
\]
and 
\[
G(n)=-\frac{\theta}{n+1}\{ \log{(n+1)}+C\}, ~~~~(1>\theta>0).
\]
But because from this it easily becomes clear, that it is 
\[
\lim{G(n)}=0, ~~~~~ [n= \infty],
\]
it is in total 
\[
\sum_{i=1}^{\infty} (-1)^{i-1} \cdot \frac{\sin{ia} \cdot \log{i}}{i}+C \cdot \sum_{i=1}^{\infty} (-1)^{i-1}  \frac{\sin{ia}}{i}= \frac{1}{2}\pi\log\bigg\{{\frac{ \Gamma({\frac{1}{2}-\frac{a}{2 \pi}})}{\Gamma({\frac{1}{2}+\frac{a}{2 \pi}})}}\bigg\}-\frac{1}{2}a\log{2\pi},
\]
whence, while
\[
\sum_{i=1}^{\infty} (-1)^{i-1} \cdot  \frac{\sin{ia}}{i}= \frac{a}{2}, ~~~~~ (a <\pi),
\]
we will finally have this remarkable formula:
\[
64.~~~ \frac{\sin{a}\cdot \log{1}}{1}-\frac{\sin{2a}\cdot \log{2}}{2}+\frac{\sin{3a}\cdot \log{3}}{3}-\frac{\sin{4a}\cdot \log{4}}{4}+\text{etc.}
\]
\[
=\frac{1}{2}\pi\log\bigg\{{\frac{ \Gamma({\frac{1}{2}-\frac{a}{2 \pi}})}{\Gamma({\frac{1}{2}+\frac{a}{2 \pi}})}}\bigg\}-\frac{1}{2}a(C+\log{2\pi});
\]
and, if $\pi-a$ is put in the place of $a$:
\[
65.~~~ \frac{\sin{a}\cdot \log{1}}{1}+\frac{\sin{2a}\cdot \log{2}}{2}+\frac{\sin{3a}\cdot \log{3}}{3}+\frac{\sin{4a}\cdot \log{4}}{4}+\text{etc.}
\]
\[
=\frac{1}{2}\pi\log\bigg\{{\frac{ \Gamma({\frac{a}{2 \pi}})}{\Gamma({\frac{1}{2}-\frac{a}{2 \pi}})}}\bigg\}-\frac{1}{2}(\pi-a)(C-\log{2\pi}).
\]
From these two formulas, by adding, we will obtain ths third one
\[
66.~~~ \frac{\sin{a}\cdot \log{1}}{1}+\frac{\sin{3a}\cdot \log{3}}{3}+\frac{\sin{5a}\cdot \log{5}}{5}+\text{etc.}
\]
\[
=\frac{1}{2}\pi\log\bigg\{{\frac{ \Gamma({\frac{a}{2 \pi}})}{\Gamma({\frac{1}{2}+\frac{a}{2 \pi}})}}\bigg\}-\frac{1}{4}\pi(C-\log{\frac{2\pi}{\tan{\frac{1}{2}a}}}).
\]
Ex. 1. After having put $a=\frac{1}{2}\pi$, it is
\[
\frac{\log{1}}{{1}}-\frac{\log{3}}{{3}}+\frac{\log{5}}{{5}}-\frac{\log{7}}{{7}}+\text{etc.}= \frac{\pi}{4}(\log{\pi}-C)-\pi \log{\Gamma{(\frac{3}{4})}}
\]
and from this
\[
\frac{{1}\cdot 5^{\frac{1}{{5}}} \cdot 9^{\frac{1}{9}} \cdot 13^{\frac{1}{{13}}} \cdot \cdots}{3^{\frac{1}{{3}}}\cdot 7^{\frac{1}{{7}}} \cdot 11^{\frac{1}{{11}}} \cdot 15^{\frac{1}{{15}}} \cdot \cdots}=\bigg\{\frac{\pi e^{-C}}{(\Gamma{(\frac{3}{4})})^4}\bigg\}^{\frac{\pi}{4}};
\]\\
we already found this above.\\

Ex. 2. After having put $a=\frac{1}{3}\pi$, it is
\begin{align*}
&\frac{\log{1}}{{1}}-\frac{\log{2}}{{2}}+\frac{\log{4}}{{4}}-\frac{\log{5}}{{5}}+\frac{\log{7}}{{7}}-\text{etc.}~~~=\frac{\pi}{\sqrt{3}}\log\bigg\{{\frac{\Gamma{(\frac{1}{3})}}{\Gamma{(\frac{2}{3})}}}\bigg\}-\frac{\pi}{3\sqrt{3}}(C+\log{2 \pi}),\\
&\frac{\log{1}}{{1}}+\frac{\log{2}}{{2}}-\frac{\log{4}}{{4}}-\frac{\log{5}}{{5}}+\frac{\log{7}}{{7}}+\text{etc.}~~~=\frac{\pi}{3\sqrt{3}}(\log{2 \pi}-2C)-\frac{2 \pi}{\sqrt{3}}\log{\Gamma{(\frac{2}{3})}},\\
&\frac{\log{1}}{{1}}-\frac{\log{5}}{{5}}+\frac{\log{7}}{{7}}-\frac{\log{11}}{{11}}+\frac{\log{13}}{{13}}-\text{etc.}=\frac{\pi}{2\sqrt{3}}\bigg\{ \log{\frac{2 \pi}{\sqrt{3}}}-C-2\log[\Gamma{(\frac{5}{6})}\cdot \Gamma{(\frac{2}{3})}]\bigg\},
\end{align*}
(of which formulas we already found the first and the third above), and from this
\begin{align*}
&\frac{{1}\cdot 4^{\frac{1}{{4}}} \cdot 7^{\frac{1}{7}} \cdot 10^{\frac{1}{{10}}} \cdot \cdots}{2^{\frac{1}{{2}}}\cdot 5^{\frac{1}{{5}}} \cdot 8^{\frac{1}{{8}}} \cdot 11^{\frac{1}{{11}}} \cdot \cdots}~~~~=\bigg\{\frac{\Gamma{(\frac{2}{3})}\cdot e^{-\frac{1}{3}C}}{\Gamma{(\frac{2}{3})}\cdot (2\pi)^{\frac{1}{3}}} \bigg\}^{\frac{\pi}{\sqrt{3}}},\\
&\frac{{1}\cdot 2^{\frac{1}{{2}}} \cdot 7^{\frac{1}{7}} \cdot 8^{\frac{1}{{8}}} \cdot \cdots}{4^{\frac{1}{{4}}}\cdot 5^{\frac{1}{{5}}} \cdot 10^{\frac{1}{{10}}} \cdot 11^{\frac{1}{{11}}} \cdot \cdots}~~=\bigg\{\frac{(2\pi)^{\frac{1}{3}}\cdot e^{-\frac{2}{3}C}}{(\Gamma{(\frac{5}{6})})^2} \bigg\}^{\frac{\pi}{\sqrt{3}}},\\
&\frac{{1}\cdot 7^{\frac{1}{{7}}} \cdot 13^{\frac{1}{13}} \cdot 19^{\frac{1}{{19}}} \cdot \cdots}{5^{\frac{1}{{5}}}\cdot 11^{\frac{1}{{11}}} \cdot 17^{\frac{1}{{17}}} \cdot 23^{\frac{1}{{23}}} \cdot \cdots}=\bigg\{\frac{(2\pi)^{\frac{1}{2}}\cdot e^{-\frac{1}{2}C}}{3^{\frac{1}{4}}\cdot \Gamma({\frac{2}{3})}\Gamma({\frac{5}{6})}} \bigg\}^{\frac{\pi}{\sqrt{3}}}.
\end{align*}
Ex. 3. After having put $a=\frac{1}{4}\pi$ in (66.), it is, of course,
\[
\frac{\log{1}}{{1}}+\frac{\log{3}}{{3}}-\frac{\log{5}}{{5}}-\frac{\log{7}}{{7}}+\text{etc.}=\frac{\pi}{\sqrt{2}}\log\bigg\{{  \frac{2^{\frac{1}{4}}\cdot \Gamma{(\frac{1}{8})}\cdot \Gamma{(\frac{3}{8})}}{e^{-\frac{1}{2}C}\cdot (2 \pi)^{\frac{3}{2}}}}\bigg\},
\]
and from this
\[
\frac{{1}\cdot 3^{\frac{1}{{3}}} \cdot 9^{\frac{1}{9}} \cdot 11^{\frac{1}{{11}}} \cdot \cdots}{5^{\frac{1}{{5}}}\cdot 7^{\frac{1}{{7}}} \cdot 13^{\frac{1}{{13}}} \cdot 15^{\frac{1}{{15}}} \cdot \cdots}= \bigg\{{  \frac{2^{\frac{1}{4}}\cdot \Gamma{(\frac{1}{8})}\cdot \Gamma{(\frac{3}{8})}}{e^{-\frac{1}{2}C}\cdot (2 \pi)^{\frac{3}{2}}}}\bigg\}^{\frac{\pi}{\sqrt{3}}}.
\]\\

\section*{~~~~~~~~~~~~~~~~~~~~~~~~~~~~~~~~~§9}

Now let us recall formula (1.) (we already proved, that it is also valid for $x=0$), which we want to multiply on both sides by
\[
\frac{e^{au}+e^{-au}}{e^{\pi u}+e^{-\pi u}}\cdot \diff{d}u;
\]
then, after having integrated between $u=0$ and $u= \infty$, it is 
\[
\int_0^{\infty} \frac{e^{au}+e^{-au}}{e^{\pi u}+e^{-\pi u}}\log{(x^2+u^2)}\diff{d}u= 2\int_0^{\infty} \frac{\diff{d}z}{z} \cdot \int_0^{\infty}\frac{e^{au}+e^{-au}}{e^{\pi u}+e^{-\pi u}} (e^{-z}-e^{-xz}\cos{uz})\diff{d}u,
\]
whence it follows by means of a known formula \footnote{See Exerc. d. Calc. Intég. par Legendre, book 2, page 186, C.J. Malmstèn} 
\[
\int_0^{\infty} \frac{e^{au}+e^{-au}}{e^{\pi u}+e^{-\pi u}}\log{(x^2+u^2)}\diff{d}u= \int_0^{\infty} \frac{\diff{d}z}{z}\bigg(\sec{\frac{1}{2}a}-\frac{2e^{-(x-\frac{1}{2})z}(1+e^{-z})\cos{\frac{1}{2}a}}{1+2e^{-z}\cos{a}+e^{-2z}}\bigg)
\]
or, after having put
\[
e^{-z}=y, ~~~\text{whence}~~~ z=\log{(\frac{1}{y})}, ~~ e^{-z}\diff{d}z= -\diff{d}y;
\]
\[
67. ~~\int_0^{\infty} \frac{e^{au}+e^{-au}}{e^{\pi u}+e^{-\pi u}}\log{(x^2+u^2)}\diff{d}u=\int_0^{1} \bigg(\sec{\frac{1}{2}a}-\frac{y^{x-1}(1+y)\cos{\frac{1}{2}a}}{1+2y\cos{a}+y^2}\bigg)\cdot \frac{\diff{d}y}{\log{\frac{1}{y}}}
\]
\[
(a<\pi)
\]
Let us call 
\[
68.~~~\mathfrak{L}(a,x)=\int_0^{\infty} \frac{e^{au}+e^{-au}}{e^{\pi u}+e^{-\pi u}}\log{(x^2+u^2)}\diff{d}u
\]
If in formula (2.) we substitute $x+\frac{1}{2}$ and $x-\frac{1}{2}$ in the place of $x$, we will have, by adding, this simple correlation between $L(a,x)$ and $\mathfrak{L}(a,x)$
\[
69. ~~~ 2\sin{\frac{1}{2}a} \cdot \mathfrak{L}(a,x)= L(a, x+\frac{1}{2})+L(a, x-\frac{1}{2}).
\]
After having perceived this, if $a$ is in any arbitrary commensurable ratio to $\pi$, i.e. for $a=\frac{m \pi}{n}$, we can easily find formulas for the function $\mathfrak{L}(a,x)$, that are analogous to those, that we proposed for $L(a,x)$ above. And for if one puts $x+\frac{1}{2}$ and $x-\frac{1}{2}$ in the place of $x$ in (7.), we will have from this
\[
L(a, x+\frac{1}{2})+L(a, x-\frac{1}{2})
\]
\[
=2\tan{\frac{1}{2}a}\cdot\log{2n}+2\sum_{i=1}^{n-1}(-1)^{i-1}\sin{ia}\log\bigg\{{\frac{\Gamma{(\frac{x+n+i+\frac{1}{2}}{2n}})}{\Gamma{(\frac{x+i+\frac{1}{2}}{2n}})}}\bigg\}
\]
\[
+2\sum_{i=1}^{n}(-1)^{i-1}\sin{ia}\log\bigg\{{\frac{\Gamma{(\frac{x+n+i-\frac{1}{2}}{2n}})}{\Gamma{(\frac{x+i-\frac{1}{2}}{2n}})}}\bigg\}
\]
\[
=2\tan{\frac{1}{2}a}\cdot\log{2n}+2\sum_{i=1}^{n-1}(-1)^{i-1}(\sin{ia}-\sin{(i-1)a})\log\bigg\{{\frac{\Gamma{(\frac{x+n+i-\frac{1}{2}}{2n}})}{\Gamma{(\frac{x+i-\frac{1}{2}}{2n}})}}\bigg\}
\]
\[
(m+n = \text{odd number})
\]
\[
L(a, x+\frac{1}{2})+L(a, x-\frac{1}{2})
\]
\[
=2\tan{\frac{1}{2}a}\cdot\log{n}+2\sum_{i=1}^{\frac{1}{2}(n-1)}(-1)^{i-1}\sin{ia}\log\bigg\{{\frac{\Gamma{(\frac{x+n-i+\frac{1}{2}}{n}})}{\Gamma{(\frac{x-i-\frac{1}{2}}{n}})}}\bigg\}
\]
\[
+2\sum_{i=1}^{\frac{1}{2}(n-1)}(-1)^{i-1}\sin{ia}\log\bigg\{{\frac{\Gamma{(\frac{x+n-i-\frac{1}{2}}{n}})}{\Gamma{(\frac{x+i+\frac{1}{2}}{n}})}}\bigg\}
\]
\[
=2\tan{\frac{1}{2}a}\cdot\log{n}+2\sum_{i=1}^{\frac{1}{2}(n-1)}(-1)^{i-1}(\sin{ia}-\sin{(i-1)a})\log\bigg\{{\frac{\Gamma{(\frac{x+n+i-\frac{1}{2}}{n}})}{\Gamma{(\frac{x+i+\frac{1}{2}}{n}})}}\bigg\}
\]
\[
(m+n = \text{even number}),
\]
because it is
\[
\sum_{i=1}^{\frac{1}{2}(n-1)}(-1)^{i-1}\sin{ia}\log\bigg\{{\frac{\Gamma{(\frac{x+n-i-\frac{1}{2}}{n}})}{\Gamma{(\frac{x+i+\frac{1}{2}}{n}})}}\bigg\}=
-\sum_{i=1}^{\frac{1}{2}(n-1)}(-1)^{i-1}\sin{(i-1)a}\log\bigg\{{\frac{\Gamma{(\frac{x+n-i-\frac{1}{2}}{n}})}{\Gamma{(\frac{x+i+\frac{1}{2}}{n}})}}\bigg\}.
\]

But hence, after having devided by $\sin{\frac{1}{2}a}$, because it is in general
\[
\sin{ia}-\sin{(i-1)a}=2\sin{\frac{1}{2}a}\cos{(i-\frac{1}{2})a},
\]
 one concludes by means of formula (69.) \\
 
70. 
$\begin{cases}
\mathfrak{L}(a,x)= \sec{\frac{1}{2}a}\log{2n}+ 2\sum_{i=1}^{n-1}(-1)^{i-1}\cos{(i-\frac{1}{2})a}\log\bigg\{{\frac{\Gamma{(\frac{x+n+i-\frac{1}{2}}{2n}})}{\Gamma{(\frac{x+i-\frac{1}{2}}{2n}})}}\bigg\}\\[5mm]
~~~~~~~~~~~(m+n= \text{odd number})\\[5mm]
\mathfrak{L}(a,x)= \sec{\frac{1}{2}a}\log{n}+\sum_{i=1}^{\frac{1}{2}(n-1)}(-1)^{i-1}\cos{(i-\frac{1}{2})a}\log\bigg\{{\frac{\Gamma{(\frac{x+n+i-\frac{1}{2}}{n}})}{\Gamma{(\frac{x+i+\frac{1}{2}}{n}})}}\bigg\}\\[5mm]
~~~~~~~~~~~(m+n= \text{even number})
\end{cases}$\\[5mm]

\section*{~~~~~~~~~~~~~~~~~~~~~~~~~~~~~~~~~§10}

If in the formulas (70.), for $x=0$, we put
\[
e^{\frac{\pi u}{n}}= y, ~~~ \text{whence} ~~~ u= \frac{n}{\pi}\log{y},
\]
of course, we will have, while $a=\frac{m \pi}{m}$,
\[
\int_1^{\infty} \frac{y^{m-1}+y^{-m-1}}{y^n+y^{-n}}\bigg\{ \log{\frac{n}{\pi} + \log{(\log{y})}}\bigg\}\diff{d}y
\]
\[
=\frac{\pi}{2n}\sec{\frac{m \pi}{2n}}\log{2n}+ \frac{\pi}{n} \cdot \sum_{i=1}^{n-1}(-1)^{i-1}\cos{(i-\frac{1}{2})\frac{m \pi}{n}} \cdot \log\bigg\{{\frac{\Gamma({\frac{1}{2}+\frac{i-\frac{1}{2}}{2n}})}{\Gamma{(\frac{i-\frac{1}{2}}{2n}})}}\bigg\}
\]
\[
(m+n= \text{odd number})
\]
\[
\int_1^{\infty} \frac{y^{m-1}+y^{-m-1}}{y^n+y^{-n}}\bigg\{ \log{\frac{n}{\pi} + \log{(\log{y})}}\bigg\}\diff{d}y
\]
\[
=\frac{\pi}{2n}\sec{\frac{m \pi}{2n}}\log{n}+ \frac{\pi}{n} \cdot \sum_{i=1}^{\frac{1}{2}(n-1)}(-1)^{i-1}\cos{(i-\frac{1}{2})\frac{m \pi}{n}} \cdot \log\bigg\{{\frac{\Gamma({1-\frac{i-\frac{1}{2}}{n}})}{\Gamma{(\frac{i-\frac{1}{2}}{n}})}}\bigg\}
\]
\[
(m+n= \text{even number})
\]
Hence indeed, because it is
\[
\int_1^{\infty} \frac{y^{m-1}+y^{-m-1}}{y^n+y^{-n}}\diff{d}y= \frac{\pi}{2n}\sec{\frac{m \pi}{2n}},
\]
it follows\\

71.
$\begin{cases}
~~~~~~~~~~~~~~~~~~~~~\int_1^{\infty} \frac{y^{m-1}+y^{-m-1}}{y^n+y^{-n}}\cdot \log{(\log{y})} \diff{d}y\\[5mm]
=\frac{\pi}{2n}\sec{\frac{m \pi}{2n}}\log{2\pi}+ \frac{\pi}{n} \cdot 2\sum_{i=1}^{n-1}(-1)^{i-1}\cos{(i-\frac{1}{2})\frac{m \pi}{n}} \cdot \log\bigg\{{\frac{\Gamma({\frac{1}{2}+\frac{i-\frac{1}{2}}{2n}})}{\Gamma{(\frac{i-\frac{1}{2}}{2n}})}}\bigg\}\\[5mm]
~~~~~~~~~~~~~~~~~~~~~(m+n= \text{odd number})\\[5mm]
~~~~~~~~~~~~~~~~~~~~~\int_1^{\infty} \frac{y^{m-1}+y^{-m-1}}{y^n+y^{-n}}\cdot \log{(\log{y})} \diff{d}y\\[5mm]
=\frac{\pi}{2n}\sec{\frac{m \pi}{2n}}\log{\pi}+ \frac{\pi}{n} \cdot \sum_{i=1}^{\frac{1}{2}(n-1)}(-1)^{i-1}\cos{(i-\frac{1}{2})\frac{m \pi}{n}} \cdot \log\bigg\{{\frac{\Gamma({1-\frac{i-\frac{1}{2}}{n}})}{\Gamma{(\frac{i-\frac{1}{2}}{n}})}}\bigg\}\\[5mm]
~~~~~~~~~~~~~~~~~~~~~(m+n= \text{even number}).
\end{cases}$\\[5mm]

After having put $m=1$ and $n=3$ in the second, and having changed $y^2$ into $y$, we will have after a certain very easy reduction
\[
72. ~~~ \int_1^{\infty} \frac{\log{(\log{y})} \diff{d}y}{1-y+y^2}= \frac{2 \pi}{\sqrt{3}}\bigg( \frac{5}{6} \log{2 \pi} -\log{\Gamma{(\frac{1}{6}})}\bigg)
\]

\section*{~~~~~~~~~~~~~~~~~~~~~~~~~~~~~~~~~§11}

In the assumption $a=\frac{m \pi}{n}$, while $m+n$ is an odd integer number, it is known
\[
\tan{\frac{1}{2}a}= \sum_{i=1}^n (-1)^{i-1} \sin{ia},
\]
and, after having put $i-1$ in the place of $i$, 
\[
\tan{\frac{1}{2}a}= -\sum_{i=1}^n (-1)^{i-1} \sin{(i-1)a};
\]
hence by adding, after having divided by $2\sec{\frac{1}{2}a}$, it is in total
\[
\sec{\frac{1}{2}a}= -\sum_{i=1}^n (-1)^{i-1} \cos{(i-\frac{1}{2})a}.
\]
Therefore the first of the formulas (70.) can also be exhibited like this:
\[
74. ~~ \mathfrak{L}(a,x) =2 \sum_{i=1}^{n}(-1)^{i-1}\cos{(i-\frac{1}{2})a} \cdot \log\bigg\{{\frac{(2n)^{\frac{1}{2}} \Gamma({\frac{1}{2}+\frac{x+i-\frac{1}{2}}{2n}})}{\Gamma{(\frac{x+i-\frac{1}{2}}{2n}})}}\bigg\}
\]
\[
(m+n= \text{odd number}).
\]
Put $x+n$ in the place of $x$ here; it is by adding, of course
\[
\mathfrak{L}(a,x+n)+\mathfrak{L}(a,x)= 2 \sum_{i=1}^{n}(-1)^{i-1}\cos{(i-\frac{1}{2})a}  \log{(x+i-\frac{1}{2})}
\]
\[
(m+n= \text{odd number}).
\]
Moreover from the same formula, if $n-x$ is substituted in the place of $x$ it emerges by adding:
\[
\mathfrak{L}(a,x)+\mathfrak{L}(a,n-x)= 2 \sum_{i=1}^{n}(-1)^{i-1}\cos{(i-\frac{1}{2})a} \cdot \log\bigg\{{\frac{(2n)^{\frac{1}{2}} \cdot \Gamma({\frac{1}{2}+\frac{x+i-\frac{1}{2}}{2n}})}{\Gamma{(\frac{x+i-\frac{1}{2}}{2n}})}}\bigg\}
\]
\[
+2 \sum_{i=1}^{n}(-1)^{i-1}\cos{(i-\frac{1}{2})a} \cdot \log\bigg\{{\frac{(2n)^{\frac{1}{2}} \cdot \Gamma({\frac{1}{2}+\frac{n+i-x-\frac{1}{2}}{2n}})}{\Gamma{(\frac{n+i-x-\frac{1}{2}}{2n}})}}\bigg\},
\]
or, because it is 
\[
\sum_{i=1}^{n}(-1)^{i-1}\cos{(i-\frac{1}{2})a} \cdot \log\bigg\{{\frac{(2n)^{\frac{1}{2}} \cdot \Gamma({\frac{1}{2}+\frac{n+i-x-\frac{1}{2}}{2n}})}{\Gamma{(\frac{n+i-x-\frac{1}{2}}{2n}})}}\bigg\}
\]
\[
=\sum_{i=1}^{n}(-1)^{i-1}\cos{(i-\frac{1}{2})a} \cdot \log\bigg\{{\frac{(2n)^{\frac{1}{2}} \cdot \Gamma({\frac{3}{2}-\frac{x+i-\frac{1}{2}}{2n}})}{\Gamma{(1-\frac{x+i-\frac{1}{2}}{2n}})}}\bigg\},
\]
also, having regard to he known property of the function $\Gamma$ (in the same way as above to the deduction of formula (22.)):
\[
76.~~ \mathfrak{L}(a,x)+\mathfrak{L}(a,n-x)= 2 \sum_{i=1}^{n-1}(-1)^{i-1}\cos{(i-\frac{1}{2})a}  \log\bigg[{(i-x-\frac{1}{2})}\cdot \cot{\frac{(i-1-\frac{1}{2})\pi}{2n}}\bigg]
\]
\[
(m+n= \text{odd number}).
\]
From the second of the formulas (70.) we can also derive analogous relations. And for, having substituted $x+n$ for $x$ there, by subtracting it will be
\[
77. ~~ \mathfrak{L}(a,x+n)-\mathfrak{L}(a,x)= 2 \sum_{i=1}^{\frac{1}{2}(n-1)}(-1)^{i-1}\cos{(i-\frac{1}{2})a}  \log\bigg\{{\frac{x+n+\frac{1}{2}-i}{x+i-\frac{1}{2}}}\bigg\}
\]
\[
(m+n= \text{even number}).
\]
And in in the same formula $n-x$ is put in the place of $x$, it also is by subtracting, after havig done the calculation,
\[
78. ~~ \mathfrak{L}(a,x)-\mathfrak{L}(a,n-x)= 2 \sum_{i=1}^{\frac{1}{2}(n-1)}(-1)^{i-1}\cos{(i-\frac{1}{2})a} \cdot \log\bigg\{{\frac{x+n+\frac{1}{2}-i}{x+i-\frac{1}{2}}\cdot \frac{\sin{\frac{(x+i-\frac{1}{2})\pi}{n}}}{\sin{\frac{(x-i+\frac{1}{2})\pi}{n}}}}\bigg\}
\]
\[
(m+n= \text{even number}).
\]
From the formulas (75. and 77.), combined with (76. and 78.), it follows, that, whether $m+n$ is odd or not, we know the function $\mathfrak{L}(a,x)$ for each value of $x$, if it is only known throughout the whole period from $x=0$ to $x=\frac{1}{2}n$. If in (76.) one puts $x=\frac{1}{2}n$, it emerges
\[
79. ~~~~ \mathfrak{L}(a,\frac{1}{2}n)= \sum_{i=1}^{n}(-1)^{i-1}\cos{(i-\frac{1}{2})a} \log\bigg\{ (\frac{1}{2}(n+1)\cdot \cot{(\frac{\pi}{4}-\frac{(i-\frac{1}{2}) \cdot \pi}{2n})}\bigg\}
\]
\[
(m+n= \text{odd number}),
\]
whence, if you at the same time consider (75.), one concludes, while $m+n$ is an odd number, that $\mathfrak{L}(a,x)$ for $x=\frac{1}{2}(2i+1)n$ can always be expressed by means of logarithms and circular functions. This formula (79.) for $i=\frac{1}{2}$ presents the expression $\log({0 \cdot \infty})$, but its true value is easily found to be $\log{(\frac{2n}{\pi})}$. \\

Now let us successively put in the place of $x$: 
\[
x,~~~x+\frac{2n}{r}, ~~~x+\frac{4n}{r}, ~~~x+\frac{6n}{r},\cdots, ~~~x+\frac{(r-1)\cdot2n}{r}
\]
in the first of the formulas (70.), and 
\[
x,~~~x+\frac{n}{r}, ~~~x+\frac{2n}{r}, ~~~x+\frac{3n}{r},\cdots, ~~~x+\frac{(r-1)\cdot n}{r}
\]
in the second; in the same way as above in §.3 from the known property of the function $\Gamma$ we will obtain\\

80.
$\begin{cases}
~~~\mathfrak{L}(a,x)+\mathfrak{L}(a, x+\frac{2n}{r})+\mathfrak{L}(a, x+\frac{4n}{r})+\cdots+\mathfrak{L}(a, x+\frac{(r-1)\cdot 2n}{r})\\[5mm]
=r\sec{\frac{1}{2}a}\cdot\log{2n}+2\sum_{i=1}^{n}(-1)^{i-1}\cos{(i-\frac{1}{2})a}\log\bigg\{{\frac{\Gamma({\frac{1}{2}r+\frac{r(x+i-\frac{1}{2})}{2n}})}{r^{\frac{1}{2}r}\Gamma{(\frac{r(x+i-\frac{1}{2})}{2n}})}}\bigg\}\\[5mm]
~~~~~~~~~~~~~~~~~~(m+n= \text{odd number})\\[5mm]
~~~\mathfrak{L}(a,x)+\mathfrak{L}(a, x+\frac{n}{r})+\mathfrak{L}(a, x+\frac{2n}{r})+\cdots+\mathfrak{L}(a, x+\frac{(r-1)\cdot n}{r})\\[5mm]
=r\sec{\frac{1}{2}a}\cdot\log{n}+2\sum_{i=1}^{\frac{1}{2}(n-1)}(-1)^{i-1}\cos{(i-\frac{1}{2})a}\log\bigg\{{\frac{\Gamma({r+\frac{r(x-i+\frac{1}{2})}{n}})}{r^{r-\frac{r(2i-1)}{n}}\cdot \Gamma{(\frac{r(x+i-\frac{1}{2})}{n}})}}\bigg\}\\[5mm]
~~~~~~~~~~~~~~~~~~(m+n= \text{even number})\\
\end{cases}$\\[2mm]

whece it becomes clear, that the first of these sums, while $r$ is an even number, can always be expressed in finite terms by means of logarithms.

\section*{~~~~~~~~~~~~~~~~~~~~~~~~~~~~~~~~~§12}

If we multiply formula (27.) on both sides by
\[
\frac{e^{au}+e^{-au}}{e^{\pi u}+e^{-\pi u}}\cdot \diff{d}u,
\]
then, after having integrated between $u=0$ and $u= \infty$, it is 
\[
\int_0^{\infty} \frac{e^{au}+e^{-au}}{e^{\pi u}+e^{-\pi u}}\cdot \frac{\cos{(s\cdot \arctan{\frac{u}{x}})}\diff{d}u}{(x^2+u^2)^{\frac{1}{2}s}}
\]
\[
=\frac{1}{\Gamma{(s)}}\cdot \int_0^{\infty} e^{-xz} \cdot z^{s-1} \cdot \int_0^{\infty} \frac{e^{au}+e^{-au}}{e^{\pi u}+e^{-\pi u}} \cos{uz} \cdot \diff{d}u
\]
\[
=\frac{\cos{\frac{1}{2}a}}{\Gamma{(s)}} \cdot \int_0^{\infty} \frac{z^{s-1} \cdot e^{-(x-\frac{1}{2})}(1+e^{-z})\cdot e^{-z} \diff{d}z}{1+2e^{-z} \cos{a}+e^{-2z}},
\]
whence, if we put $e^{-z}=y$ on the right hand side, it emerges
\[
81. ~~~ \int_0^{\infty} \frac{e^{au}+e^{-au}}{e^{\pi u}+e^{-\pi u}}\cdot \frac{\cos{(s\cdot \arctan{\frac{u}{x}})}\diff{d}u}{(x^2+u^2)^{\frac{1}{2}s}} = \frac{\cos{\frac{1}{2}a}}{\Gamma{(s)}} \cdot \int_0^1 \frac{y^{x-\frac{1}{2}}(1+y) (\log{\frac{1}{y}})^{s-1}\diff{d}y}{1+2y\cos{a}+y^2} 
\]
and for $x=0$ (if only $1>s>0$), after having changed $y$ into $y^2$, 
\[
82. ~~~ \int_0^{\infty} \frac{e^{au}+e^{-au}}{e^{\pi u}+e^{-\pi u}}\cdot \frac{\diff{d}u}{u^s}= \frac{2^s \cos{\frac{1}{2}a}}{\cos{\frac{1}{2}s \pi} \Gamma{(s)}} \cdot \int_0^1 \frac{(1+y^2) (\log{\frac{1}{y}})^{s-1}\diff{d}y}{1+2y^2\cos{a}+y^4}. 
\]
But from this formula, if $a=\frac{m \pi}{n}$ and one puts $e^{-\frac{\pi u}{n}}=y$, after having done the transformation, it is
\[
83. ~~~ \int_0^1 \frac{y^{m-1}+y^{-m-1}}{y^n+y^{-n}} \cdot \frac{\diff{d}y}{(\log{\frac{1}{y}})^s}= \frac{2 \cdot (\frac{\pi}{n})^{1-s} \cos{\frac{m \pi}{2n}}}{\cos{\frac{1}{2}s \pi} \cdot \Gamma{(s)}} \cdot  \int_0^1 \frac{(1+y^2) (\log{\frac{1}{y}})^{s-1}\diff{d}y}{1+2y^2\cos{\frac{m \pi}{n}}+y^4}
\]\\

After having put $m=1$ and $n=2$, if we call for the sake of brevity
\[
84. ~~~ Q(s)= \int_0^1 \frac{1+y^2}{1+y^4} \cdot (\log{\frac{1}{y}})^{s-1}\diff{d}y,
\]
we will have this relation between the function $Q(s)$ and its compliment:
\[
85. ~~~ Q(1-s)= \frac{2 \cdot (\frac{\pi}{4})^{1-s} \cdot \sin{\frac{ \pi}{4}}}{\cos{\frac{1}{2}s \pi} \cdot \Gamma{(s)}} \cdot Q(s)
\]
After having further put $m=1$ and $n=3$ in (83.) and having done the reductions, it emerges
\[
86. ~~ \int_0^1 \frac{(\log{\frac{1}{y}})^{-s}\diff{d}y}{1-y+y^2}=\frac{ (\frac{2}{3}\pi)^{1-s} \cdot \sin{\frac{ 1}{3}\pi}}{\cos{\frac{1}{2}s \pi} \cdot \Gamma{(s)}} \cdot \int_0^1 \frac{2(1+y^2)}{1+y^2+y^4} \cdot (\log{\frac{1}{y}})^{s-1}\diff{d}y.
\]
But because it is
\[
\int_0^1 \frac{2(1+y^2)}{1+y^2+y^4} \cdot (\log{\frac{1}{y}})^{s-1}\diff{d}y=\int_0^1 \frac{(\log{\frac{1}{y}})^{s-1}\diff{d}y}{1+y+y^2} + \int_0^1 \frac{(\log{\frac{1}{y}})^{s-1}\diff{d}y}{1-y+y^2},
\]
and also
\[
\int_0^1 \frac{(\log{\frac{1}{y}})^{s-1}\diff{d}y}{1+y+y^2}= \int_0^1 \frac{2^{2s} \cdot y^3(\log{\frac{1}{y}})^{s-1}\diff{d}y}{1+y^4+y^8}
\]
\[
= 2^{2s-2}\int_0^1 \frac{(\log{\frac{1}{y}})^{s-1}\diff{d}y}{1+y+y^2} -2^{2s-2}\int_0^1 \frac{2y(\log{\frac{1}{y}})^{s-1}\diff{d}y}{1+y^2+y^4}
\]
i.e.
\[
(1+2^{s-1})\int_0^1 \frac{(\log{\frac{1}{y}})^{s-1}\diff{d}y}{1+y+y^2}=2^{s-1} \cdot \int_0^1 \frac{(\log{\frac{1}{y}})^{s-1}\diff{d}y}{1-y+y^2};
\]
it also is from (87.)
\[
\int_0^1 \frac{2(1+y^2)}{1+y+y^2} \cdot (\log{\frac{1}{y}})^{s-1}\diff{d}y = \frac{1+2^s}{1+2^{s-1}} \cdot \int_0^1 \frac{(\log{\frac{1}{y}})^{s-1}\diff{d}y}{1-y+y^2},
\]
which substituted in (86.), after having put for the sake of brevity
\[
88. ~~~~ R(s)= \int_0^1 \frac{(\log{\frac{1}{y}})^{s-1}\diff{d}y}{1-y+y^2},
\]
gives this remarkable relation:
\[
89. ~~~~ R(1-s)= \frac{ (\frac{1}{3}\pi)^{1-s} \cdot \sin{\frac{ 1}{3}\pi}}{\cos{\frac{1}{2}s \pi} \cdot \Gamma{(s)}} \cdot \frac{1+2^s}{1+2^{s-1}} \cdot R(s).
\]
From the formulas (85. and 89.), by taking logarithms, we will obtain
\[
\log{Q(1-s)}= \log{2}\sin{\frac{\pi}{4}}+(1-s) \log{\frac{\pi}{4}}- \log{\cos{\frac{1}{2}s \pi}}- \log{\Gamma{(s)}}+\log{Q(s)},
\]
\[
\log{R(1-s)}= (1-s) \log{\frac{1}{3}\pi}+\log{\sin{\frac{1}{3}\pi}}+\log{(1+2^s)}-\log{\cos{\frac{1}{2}s \pi}}
\]
\[
-\log{\Gamma{(s)}}-\log{(1+2^{s-1})}+\log{Q(s)},
\]
whence, after having put for the sake of brevity
\begin{align*}
M(s)&=\frac{\diff{d} Q(s)}{\diff{d}s}= \int_0^1 \frac{1+y^2}{1+y^4} \cdot \log{(\log\frac{1}{y})} \cdot (\log{\frac{1}{y}})^{s-1}\diff{d}y\\
N(s)&= \frac{\diff{d} R(s)}{\diff{d}s} = \int_0^1 \frac{\log{(\log\frac{1}{y})}\cdot (\log{\frac{1}{y}})^{s-1}\diff{d}y }{1-y+y^2},
\end{align*}
by differentiating, we will have\\[5mm]

90.
$\begin{cases}
\frac{M(s)}{Q(s)}+\frac{M(1-s)}{Q(1-s)}= \log{\frac{\pi}{4}}- \frac{1}{2}\pi \tan{\frac{1}{2}s \pi}+Z^{\prime}(s) ~~~~~~ \text{and}\\[5mm]
\frac{N(s)}{R(s)}+\frac{N(1-s)}{R(1-s)}~= \log{\frac{\pi}{3}}- \frac{1}{2}\pi \tan{\frac{1}{2}s \pi}+Z^{\prime}(s)-\frac{\log{2}}{(1+2^s)(1+2^{s-1})}
\end{cases} $\\[2mm]

Let us suppose $s=\frac{1}{2}$, then it will be

\[
\int_0^1 \frac{1+y^2}{1+y^4} \cdot \frac{\log{(\log\frac{1}{y})} \diff{d}y}{\sqrt{\log{\frac{1}{y}}}}= \frac{1}{2}\bigg(\log{\frac{\pi}{16}}-\frac{1}{2}\pi-C\bigg) \cdot \int_0^1 \frac{1+y^2}{1+y^4} \cdot \frac{ \diff{d}y}{\sqrt{\log{\frac{1}{y}}}}
\]
\[
\int_0^1 \frac{1}{1-y+y^2}\cdot \frac{\log{(\log\frac{1}{y})}\diff{d}y}{\sqrt{\log{\frac{1}{y}}}}= \frac{1}{2}\bigg(\log{\frac{1}{3}\pi}-\frac{1}{2}\pi-C-(5-2\sqrt{2})\log{2}\bigg) \cdot \int_0^1 \frac{1}{1-y+y^2}\cdot \frac{\diff{d}y}{\sqrt{\log{\frac{1}{y}}}},
\]
or, by putting $\log{\frac{1}{y}}=x$:\\[5mm]
91.
$\begin{cases}
\int_0^{\infty} \frac{e^x+e^{-x}}{e^{2x}+e^{-2x}} \cdot \frac{\log{x} \cdot \diff{d}x}{\sqrt{x}}= \frac{1}{2}\bigg(\log{\frac{\pi}{16}}-\frac{1}{2}\pi-C\bigg) \cdot \int_0^{\infty} \frac{e^x+e^{-x}}{e^{2x}+e^{-2x}} \cdot \frac{ \diff{d}x}{\sqrt{x}},\\[5mm]
\int_0^{\infty} \frac{\log{x}}{e^x-1+e^{-x}} \cdot \frac{\diff{d}x}{\sqrt{x}}= \frac{1}{2}\bigg(\log{\frac{1}{3}\pi}-\frac{1}{2}\pi-C-(5-2\sqrt{2})\log{2}\bigg) \cdot \int_0^{\infty} \frac{\diff{d}x}{e^x-1+e^{-x}} \cdot \frac{1}{\sqrt{x}},
\end{cases}$\\[5mm]

which formulas are of completely the same kind as the formulas (38.).\\

Form these formulas by exactly the same method, that we used to find the formulas (42.), it is possible to deduce these equally remarkable relations:\\[5mm]

92. 
$\begin{cases}
\frac{\log{1}}{\sqrt{1}}+\frac{\log{3}}{\sqrt{3}}-\frac{\log{5}}{\sqrt{5}}-\frac{\log{7}}{\sqrt{7}}+\frac{\log{9}}{\sqrt{9}}+\text{etc.}\\[5mm]
\frac{1}{2}(\frac{1}{2}\pi-C-\log{\pi})\{\frac{1}{\sqrt{1}}+\frac{1}{\sqrt{3}}-\frac{1}{\sqrt{5}}-\frac{1}{\sqrt{7}}+\frac{1}{\sqrt{9}}-\text{etc.} \}\\[5mm]
\frac{\log{1}}{\sqrt{1}}+\frac{\log{2}}{\sqrt{2}}-\frac{\log{4}}{\sqrt{4}}-\frac{\log{5}}{\sqrt{5}}+\frac{\log{7}}{\sqrt{7}}+\frac{\log{8}}{\sqrt{8}}-\frac{\log{10}}{\sqrt{10}}-\text{etc.}\\[5mm]
\frac{1}{2}(\frac{1}{2}\pi-C-2\sqrt{2}\log{2}-\log{\frac{\pi}{6}})\{\frac{1}{\sqrt{1}}+\frac{1}{\sqrt{2}}-\frac{1}{\sqrt{4}}-\frac{1}{\sqrt{5}}+\frac{1}{\sqrt{7}}+\frac{1}{\sqrt{8}}-\text{etc.} \}
\end{cases}$\\

\section*{~~~~~~~~~~~~~~~~~~~~~~~~~~~~~~~~~§13}

If we expand 
\[
\frac{e^{au}+e^{-au}}{e^{\pi u}+e^{-\pi u}}
\]
into a series, it is identical manner, of course
\[
\frac{e^{au}+e^{-au}}{e^{\pi u}+e^{-\pi u}}= \sum_{i=0}^{n-1}(-1)^i\bigg[ e^{-[(2i+1)\pi -a]u}+e^{-[(2i+1)\pi +a]u}\bigg]+\frac{(-1)^n e^{-2n \pi u}(e^{au}+e^{-au})}{e^{\pi u}+e^{-\pi u}},
\]
whence
\[
\int_0^{\infty} \frac{e^{au}+e^{-au}}{e^{\pi u}+e^{-\pi u}} \cdot \frac{\diff{d}u}{u^s}
\]
\[
=\Gamma{(1-s)} \cdot \sum_{i=0}^{n-1}(-1)^i\bigg[ \frac{1}{((2i+1)\pi-a})^{1-s}+ \frac{1}{((2i+1)\pi+a)^{1-s}}\bigg]+(-1)^n \cdot \varphi_1{(n)},
\]
where we put for the sake of brevity 
\[
\varphi_1{(n)}=\int_0^{\infty} \frac{e^{au}+e^{-au}}{e^{\pi u}+e^{-\pi u}} \cdot \frac{e^{-2n \pi u}\diff{d}u}{u^s}= E \cdot \int_0^{\infty} \frac{e^{-2n \pi u}\diff{d}u}{u^s} = \frac{M \Gamma{(1-s)}}{(2n \pi)^{1-s}},
\]
(while $E$ is a certain finite quantity); and because it obviously is
\[
\lim{\varphi{(n)}}=0, ~~~~~~~ [n= \infty]
\]
also
\[
93. ~~~~ \int_0^{\infty} \frac{e^{au}+e^{-au}}{e^{\pi u}+e^{-\pi u}} \cdot \frac{\diff{d}u}{u^s}= \Gamma{(1-s)} \cdot \sum_{i=0}^{\infty}(-1)^i\bigg[ \frac{1}{((2i+1)\pi-a})^{1-s}+ \frac{1}{((2i+1)\pi+a)^{1-s}}\bigg].
\]
Now we will indeed, while
\[
\frac{(1+y^2)\cos{\frac{1}{2}a}}{1+y^2\cos{a}+y^4}=\sum_{i=0}^{n-1} (-1)^{i} \cdot y^{2i} \cdot \cos{(i+\frac{1}{2})a} +(-1)^n \cdot y^{2n} \cdot \frac{\cos{(n+\frac{1}{2})a}+y^2\cos{(n-\frac{1}{2})a}}{1+y^2\cos{a}+y^4},
\]
also have
\[
\int_0^1 \frac{(1+y^2)\cos{\frac{1}{2}a}}{1+y^2\cos{a}+y^4} \cdot \frac{\diff{d}y}{(\log{\frac{1}{y}})^{1-s}}
\]
\[
= \sum_{i=0}^{n-1} (-1)^{i} \cdot \frac{\cos{(i+\frac{1}{2})a}}{(2i+1)^s}+(-1)^n\bigg(p(n)\cos{(n+\frac{1}{2})a}+p(n+1)\cos{(n-\frac{1}{2})a}\bigg),
\]
after having put for the sake of brevity
\[
p(n)= \int_0^1 \frac{y^{2n} \cdot (\log{\frac{1}{y}})^{s} \diff{d}y }{1+y^2\cos{a}+y^4}=\theta \cdot \int_0^1 y^{2n}  (\log{\frac{1}{y}})^{s-1} \diff{d}y= \frac{\theta \cdot \Gamma{(s)}}{(2n+1)^s} ~~ (1> \theta > 0).
\]
Therefore it is obviously
\[
\lim{p(n)}=0 ~~~~~~ [n=8],
\]
whence it is, of course,
\[
94. ~~~~ \int_0^1 \frac{(1+y^2)\cos{\frac{1}{2}a}}{1+y^2\cos{a}+y^4} \cdot \frac{\diff{d}y}{(\log{\frac{1}{y}})^{1-s}}= \Gamma{(s)} \cdot \sum_{i=0}^{\infty} (-1)^{i} \cdot \frac{\cos{(i+\frac{1}{2})a}}{(2i+1)^s}.
\]
But after having substituted the values, that the formulas (93. and 94.) give, in (82.), if $s$ is changed into $1-s$, it emerges\\[5mm]
95.
$\begin{cases}
\frac{1}{(\pi -a)^s}+\frac{1}{(\pi +a)^s}-\frac{1}{(3\pi -a)^s}-\frac{1}{(3\pi +a)^s}+\frac{1}{(5\pi -a)^s}+\frac{1}{(5\pi +a)^s}-\text{etc.}\\[3mm]
=\frac{2^{1-s}}{\sin{\frac{1}{2}s \pi}\cdot \Gamma{(s)}}\bigg\{ \frac{\cos{\frac{1}{2}a}}{1^{1-s}}-\frac{\cos{\frac{3}{2}a}}{3^{1-s}}+\frac{\cos{\frac{5}{2}a}}{5^{1-s}}-\frac{\cos{\frac{7}{2}a}}{7^{1-s}}+\frac{\cos{\frac{9}{2}a}}{9^{1-s}}-\text{etc.}\bigg\}
\end{cases}$\\[5mm]

and if $\pi-a$ is put in the place of $a$:\\[5mm]
96.
$\begin{cases}
\frac{1}{a^s}+\frac{1}{(2\pi -a)^s}-\frac{1}{(2\pi +a)^s}-\frac{1}{(4\pi -a)^s}+\frac{1}{(4\pi +a)^s}+\frac{1}{(6\pi -a)^s}-\text{etc.}\\[3mm]
=\frac{2^{1-s}}{\sin{\frac{1}{2}s \pi}\cdot \Gamma{(s)}}\bigg\{ \frac{\sin{\frac{1}{2}a}}{1^{1-s}}+\frac{\sin{\frac{3}{2}a}}{3^{1-s}}+\frac{\sin{\frac{5}{2}a}}{5^{1-s}}+\frac{\sin{\frac{7}{2}a}}{7^{1-s}}+\frac{\sin{\frac{9}{2}a}}{9^{1-s}}+\text{etc.}\bigg\}
\end{cases}$\\[5mm]

and if we suppose $a=\frac{m \pi}{n}$ ($m<n$ integer number):\\[5mm]
97.
$\begin{cases}
\frac{1}{(n -m)^s}+\frac{1}{(n+m)^s}-\frac{1}{(3n -m)^s}-\frac{1}{(3n +m)^s}+\frac{1}{(5n -m)^s}+\frac{1}{(5n +m)^s}-\text{etc.}\\[3mm]
=\frac{2\cdot (\frac{\pi}{2n})^s}{\sin{\frac{1}{2}s \pi}\cdot \Gamma{(s)}}\bigg\{ \frac{\cos{\frac{m \pi}{2n}}}{1^{1-s}}-\frac{\cos{\frac{3m \pi}{2n}}}{3^{1-s}}+\frac{\cos{\frac{5m \pi}{2n}}}{5^{1-s}}-\frac{\cos{\frac{7m \pi}{2n}}}{7^{1-s}}+\text{etc.}\bigg\}\\[3mm]
\frac{1}{m^s}+\frac{1}{(2m -n)^s}-\frac{1}{(2m +n)^s}-\frac{1}{(4m -n)^s}+\frac{1}{(4m +n)^s}+\frac{1}{(6m -n)^s}-\text{etc.}\\[3mm]
=\frac{2\cdot (\frac{\pi}{2n})^s}{\sin{\frac{1}{2}s \pi}\cdot \Gamma{(s)}}\bigg\{ \frac{\sin{\frac{m \pi}{2n}}}{1^{1-s}}+\frac{\sin{\frac{3m \pi}{2n}}}{3^{1-s}}+\frac{\sin{\frac{5m \pi}{2n}}}{5^{1-s}}+\frac{\sin{\frac{7m \pi}{2n}}}{7^{1-s}}+\text{etc.}\bigg\}
\end{cases}$\\[5mm]

After having put $m=1$, $n=2$, it is
\[
\frac{1}{1^s}+\frac{1}{3^s}-\frac{1}{5^s}-\frac{1}{7^s}+\frac{1}{9^s}+\frac{1}{11^s}-\text{etc.}
\]
\[
=\frac{2 \cdot (\frac{\pi}{4})^s \cdot \sin{\frac{\pi}{4}}}{\sin{\frac{1}{2} s \pi} \cdot \Gamma{(s)}}\cdot\bigg\{ \frac{1}{1^{1-s}}+\frac{1}{3^{1-s}}-\frac{1}{5^{1-s}}-\frac{1}{7^{1-s}}+\text{etc.}\bigg\};
\]
we already found this above in (54.)\\

Ex. 2. After having put $m=1$, $n=3$, the first of the formulas (97.), after having done certain very easy reductions, gives the formulas (53.); but from the second, after having put for the sake of brevity
\begin{align*}
T(s)&=\frac{1}{1^s}+\frac{1}{5^s}-\frac{1}{7^s}-\frac{1}{11^s}+\frac{1}{13^s}+\frac{1}{17^s}-\text{etc.}\\
W(s)&=\frac{1}{1^s}-\frac{1}{3^s}+\frac{1}{5^s}-\frac{1}{7^s}+\frac{1}{9^s}-\frac{1}{11^s}+\text{etc.}
\end{align*}

it is 
\[
T(s)= \frac{2(\frac{\pi}{6})^s}{\sin{\frac{1}{2} s \pi} \cdot \Gamma{(s)}}\bigg\{ \frac{1}{2}T(1-s)+3^{s-1}W(1-s)\bigg\},
\]
whence, because it is
\[
W(s)=T(s) -3^{-s}W(s),
\]
we will have
\[
98. ~~ T(s)= \frac{(\frac{1}{2}\pi)^s}{\sin{\frac{1}{2} s \pi} \cdot \Gamma{(s)}}\cdot \frac{1+3^{-s}}{1+3^{s-1}} \cdot T(1-s)
\]

\section*{~~~~~~~~~~~~~~~~~~~~~~~~~~~~~~~~~§14}

Now let us differentiate formula (93.) with respect to $s$ as a variable, then it will be
\[
\int_0^{\infty} \frac{e^{au}+e^{-au}}{e^{\pi u}+e^{-\pi u}} \cdot u^{-s} \log{u} \diff{d}u
\]
\[
=Z^{\prime}(1-s) \cdot \int_0^{\infty} \frac{e^{au}+e^{-au}}{e^{\pi u}+e^{-\pi u}} \cdot \frac{\diff{d}u}{u^s} -\Gamma{(1-s)} \cdot \sum_{i=0}^{\infty} (-1)^i \bigg[\frac{\log{((2i+1)\pi-a)}}{((2i+1)\pi-a)^{1-s}}+\frac{\log{((2i+1)\pi+a)}}{((2i+1)\pi+a)^{1-s}}\bigg], 
\]  
and for $s=0$, while $Z^{\prime}(1)=-C$ and 
\[
99. ~~~~~ \int_0^{\infty} \frac{e^{au}+e^{-au}}{e^{\pi u}+e^{-\pi u}} \cdot \diff{d}u =\frac{1}{2} \sec{\frac{1}{2}a},
\]
we will have
\[
\sum_{i=0}^{\infty} (-1)^i \bigg[\frac{\log{((2i+1)\pi-a)}}{((2i+1)\pi-a)^{1-s}}+\frac{\log{((2i+1)\pi+a)}}{((2i+1)\pi+a)^{1-s}}\bigg]
\]
\[
=-\frac{1}{2}C\sec{\frac{1}{2}a}- \int_0^{\infty} \frac{e^{au}+e^{-au}}{e^{\pi u}+e^{-\pi u}} \log{u} \cdot \diff{d}u.
\]
 Whereever $a$ is indeed in a commensurable ratio to $\pi$, the integrals on the right-hand side give the value of the formula (70.) (for $x=0$). Therefore let be $a=\frac{m \pi}{n}$; then one will obtain
 \[
\frac{n}{\pi} \cdot  \sum_{i=0}^{\infty} (-1)^i \bigg[\frac{\log{((2i+1)n-m)}}{((2i+1)n-m)^{1-s}}+\frac{\log{((2i+1)n+m)}}{((2i+1)n+m)^{1-s}}\bigg]
 \] 
\[
+\log{\frac{\pi}{n}} \cdot \sum_{i=0}^{\infty} (-1)^i \bigg[\frac{1}{((2i+1)\pi-\frac{m \pi}{n})^{1-s}}+\frac{1}{((2i+1)\pi+\frac{m \pi}{n})^{1-s}}\bigg]
\]
\[
=-\frac{1}{2}C \cdot \sec{\frac{m \pi}{2n}}- \int_0^{\infty} \frac{e^{\frac{m \pi}{n} \cdot u}+e^{-\frac{m \pi}{n} \cdot u}}{e^{\pi u}+e^{-\pi u}} \cdot \log{u} \cdot \diff{d}u,
\]
whence, because it is from (93.) (for $s=0$)
\[
\sum_{i=0}^{\infty} (-1)^i \bigg[\frac{1}{(2i+1)\pi-\frac{m \pi}{n}}+\frac{1}{(2i+1)\pi+\frac{m \pi}{n}}\bigg]= \int_0^{\infty} \frac{e^{\frac{m \pi}{n} \cdot u}+e^{-\frac{m \pi}{n} \cdot u}}{e^{\pi u}+e^{-\pi u}}  \cdot \diff{d}u= \sec{\frac{m \pi}{2n}},
\]
if finally is from the cited formula (70.)\\

100.
$\begin{cases}
\frac{\log{(n-m)}}{n-m}+\frac{\log{(n+m)}}{n+m}-\frac{\log{(3n-m)}}{3n-m}-\frac{\log{(3n+m)}}{3n+m}+\frac{\log{(5n-m)}}{5n-m}+\frac{\log{(5n+m)}}{5n+m}-\text{etc.}\\[3mm]
= -\frac{\pi}{2n} \cdot \sec{\frac{m \pi}{2n}}(C+\log{2 \pi})- \frac{\pi}{n} \cdot \sum_{i=1}^{n} (-1)^{i-1} \cos{(2i+1)\frac{im \pi}{n}}\log\bigg\{{\frac{\Gamma{(\frac{1}{2}+\frac{i-\frac{1}{2}}{2n}})}{\Gamma{(\frac{i-\frac{1}{2}}{2n}})}}\bigg\}\\[5mm]
~~~~~~~~~~~~~~~~~~~~~~~~~(m+n= \text{odd number})\\[5mm]
\frac{\log{(n-m)}}{n-m}+\frac{\log{(n+m)}}{n+m}-\frac{\log{(3n-m)}}{3n-m}-\frac{\log{(3n+m)}}{3n+m}+\frac{\log{(5n-m)}}{5n-m}+\frac{\log{(5n+m)}}{5n+m}-\text{etc.}\\[3mm]
= -\frac{\pi}{2n} \cdot \sec{\frac{m \pi}{2n}}(C+\log{ \pi})- \frac{\pi}{n} \cdot \sum_{i=1}^{\frac{1}{2}(n-1)} (-1)^{i-1} \cos{(2i+1)\frac{im \pi}{n}}\log\bigg\{{\frac{\Gamma{(1-\frac{i-\frac{1}{2}}{n}})}{\Gamma{(\frac{i-\frac{1}{2}}{n}})}}\bigg\}\\[5mm]
~~~~~~~~~~~~~~~~~~~~~~~~~(m+n= \text{even number}),\\
\end{cases}$\\[5mm]

and, if $m-n$ is put in the place of $m$\\

101.
$\begin{cases}
\frac{\log{m}}{m}+\frac{\log{(2n-m)}}{2n-m}-\frac{\log{(2n+m)}}{2n+m}-\frac{\log{(4n-m)}}{4n-m}+\frac{\log{(4n+m)}}{4n+m}+\frac{\log{(6n-m)}}{6n-m}-\text{etc.}\\[3mm]
= -\frac{\pi}{2n} \cdot \csc{\frac{m \pi}{2n}}(C+\log{2 \pi})- \frac{\pi}{n} \cdot \sum_{i=1}^{n} (-1)^{i-1} \sin{(2i+1)\frac{im \pi}{n}}\log\bigg\{{\frac{\Gamma{(\frac{1}{2}+\frac{i-\frac{1}{2}}{2n}})}{\Gamma{(\frac{i-\frac{1}{2}}{2n}})}}\bigg\}\\[5mm]
~~~~~~~~~~~~~~~~~~~~~~~~~(m+n= \text{odd number})\\[5mm]
\frac{\log{m}}{m}+\frac{\log{(2n-m)}}{2n-m}-\frac{\log{(2n+m)}}{2n+m}-\frac{\log{(4n-m)}}{4n-m}+\frac{\log{(4n+m)}}{4n+m}+\frac{\log{(6n-m)}}{6n-m}-\text{etc.}\\[3mm]
= -\frac{\pi}{2n} \cdot \csc{\frac{m \pi}{2n}}(C+\log{ \pi})- \frac{\pi}{n} \cdot \sum_{i=1}^{\frac{1}{2}(n-1)} (-1)^{i-1} \sin{(2i+1)\frac{im \pi}{n}}\log\bigg\{{\frac{\Gamma{(1-\frac{i-\frac{1}{2}}{n}})}{\Gamma{(\frac{i-\frac{1}{2}}{n}})}}\bigg\}\\[5mm]
~~~~~~~~~~~~~~~~~~~~~~~~~(m+n= \text{even number}).\\
\end{cases}$\\[5mm]

Ex. 1. After having put $m=1$, $n=2$, it is
\[
\frac{\log{1}}{1}+\frac{\log{3}}{3}-\frac{\log{5}}{5}-\frac{\log{7}}{7}+\frac{\log{9}}{9}+\frac{\log{11}}{11}-\text{etc.}= \frac{\pi}{\sqrt{2}}\log\bigg\{\frac{2^{\frac{1}{4}}\cdot \Gamma{(\frac{1}{8})} \cdot \Gamma{(\frac{3}{8})}}{e^{\frac{1}{2}C}\cdot (2 \pi)^{\frac{3}{2}}}\bigg\}
\]
and from this
\[
\frac{{1}\cdot 3^{\frac{1}{{3}}} \cdot 9^{\frac{1}{9}} \cdot 11^{\frac{1}{{11}}} \cdot 17^{\frac{1}{{17}}} \cdot \cdots}{5^{\frac{1}{{5}}}\cdot 7^{\frac{1}{{7}}} \cdot 13^{\frac{1}{{13}}} \cdot 15^{\frac{1}{{25}}} \cdot  21^{\frac{1}{{21}}} \cdot \cdots}=\bigg\{\frac{2^{\frac{1}{4}}\cdot \Gamma{(\frac{1}{8})} \cdot \Gamma{(\frac{3}{8})}}{e^{\frac{1}{2}C}\cdot (2 \pi)^{\frac{3}{2}}}\bigg\}^\frac{\pi}{\sqrt{2}}.
\]
Ex. 3. $m=1$, $n=3$, it is 
\[
\frac{\log{1}}{1}+\frac{\log{2}}{2}-\frac{\log{4}}{4}-\frac{\log{5}}{5}+\frac{\log{7}}{7}+\frac{\log{8}}{8}-\text{etc.}= \frac{\pi}{3\sqrt{3}}(\log{2\pi}-2C)-\frac{2 \pi}{\sqrt{3}}\log{\Gamma{(\frac{5}{6})}},
\]
\[
\frac{\log{1}}{1}+\frac{\log{5}}{5}-\frac{\log{7}}{7}-\frac{\log{11}}{11}+\frac{\log{13}}{13}+\text{etc.}= \frac{1}{3}\pi\log\bigg\{\frac{3^{\frac{1}{4}}e^{-C} \cdot \Gamma{(\frac{1}{4})}}{2^{\frac{1}{4}}\cdot (\Gamma{(\frac{3}{4})})^3}\bigg\},
\]
and from this
\[
\frac{{1}\cdot 2^{\frac{1}{{2}}} \cdot 7^{\frac{1}{7}} \cdot 8^{\frac{1}{{8}}} \cdot \cdots}{4^{\frac{1}{{4}}}\cdot 5^{\frac{1}{{5}}} \cdot 10^{\frac{1}{{10}}} \cdot 11^{\frac{1}{{11}}}   \cdot \cdots}=\bigg\{\frac{(2\pi)^{\frac{1}{3}} \cdot e^{-\frac{2}{3}C}}{(\Gamma{(\frac{5}{6})})^2}\bigg\}^{\frac{\pi}{\sqrt{3}}},
\]
\[
\frac{{1}\cdot 5^{\frac{1}{{5}}} \cdot 13^{\frac{1}{13}} \cdot 17^{\frac{1}{{17}}} \cdot \cdots}{7^{\frac{1}{{7}}}\cdot 11^{\frac{1}{{11}}} \cdot 19^{\frac{1}{{19}}} \cdot 23^{\frac{1}{{23}}}   \cdot \cdots}=\bigg\{\frac{3^{\frac{1}{4}} \cdot e^{-C} \cdot \Gamma{(\frac{1}{4})}}{2^{\frac{1}{4}}\cdot (\Gamma{(\frac{3}{4})})^3}\bigg\}^{\frac{\pi}{3}}.
\]\\[5mm]

P.S. From the preceeding it is also possible to prove this remarkable formula
\[
\cos{\frac{1}{2}a}\cdot \log{1}-\frac{1}{3}\cos{\frac{3}{2}a}\cdot \log{3}+\frac{1}{3}\cos{\frac{5}{2}a}\cdot \log{5}-\frac{1}{7}\cos{\frac{7}{2}a}\cdot \log{7}+\frac{1}{9}\cos{\frac{9}{2}a}\cdot \log{9}-\text{etc.}
\]
\[
=\frac{\pi}{4}(\log{\pi}-C- \log{\cos{\frac{1}{2}a}})- \frac{1}{2}\pi \log\bigg[\Gamma{(\frac{3}{4}+\frac{a}{4\pi})} \cdot \Gamma{(\frac{3}{4}-\frac{a}{4\pi})}\bigg],
\]
which is valid, whatever $a< \pi$ is.\\[10mm]

Upsala 1st May 1846\\[3mm]

\end{document}